\documentclass[12pt,twoside]{amsart}
\usepackage{amssymb}

\usepackage[all]{xy}
\nonstopmode

\numberwithin{equation}{section} 

\font\tengothic=eufm10 scaled\magstep 1 \font\sevengothic=eufm7
scaled\magstep 1
\newfam\gothicfam
       \textfont\gothicfam=\tengothic
       \scriptfont\gothicfam=\sevengothic
\def\goth#1{{\fam\gothicfam #1}}

\newtheorem{theorem}{Theorem}[section]
\newtheorem{lemma}[theorem]{Lemma}
\newtheorem{proposition}[theorem]{Proposition}
\newtheorem{corollary}[theorem]{Corollary}
\newtheorem{conjecture}[theorem]{Conjecture}

\theoremstyle{definition}
\newtheorem{definition}[theorem]{Definition} 
\newtheorem{remark}[theorem]{Remark}
\newtheorem{example}[theorem]{Example}

\newtheorem{question}[theorem]{Question}

\newcommand{\codim}{\operatorname{codim}}
\newcommand{\coker}{\operatorname{coker}}
\newcommand{\aut}{\operatorname{aut}}

\newcommand{\Hom}{\operatorname{Hom}}
\newcommand{\ext}{\operatorname{ext}}
\newcommand{\Ext}{\operatorname{Ext}}

\newcommand{\depth}{\operatorname{depth}}

\newcommand{\im}{\operatorname{im}}

\newcommand{\TT}{\operatorname{Tor}}

\newcommand{\Hi}{\operatorname{Hilb}}
\newcommand{\smallbox}{\hbox{\tiny $\qed$}}

\newcommand{\proj}[1]
{ \mathchoice
            { {\mathbb P}^{#1} }
            { {\mathbb P}^{#1} }
            { {\mathbb P}^{#1} }
            { {\mathbb P}^{#1} }
          }

\newcommand{\Tor}{\operatorname{Tor}}

\newcommand{\cD}{{\mathcal D}}
\newcommand{\cA}{{\mathcal A}}
\newcommand{\cB}{{\mathcal B}}
\newcommand{\cH}{{\mathcal H}}
\newcommand{\cI}{{\mathcal I}}
\newcommand{\cK}{{\mathcal K}}

\newcommand{\cU}{{\mathcal U}}

\newcommand{\cC}{{\mathcal C}}

\newcommand{\cF}{{\mathcal F}}
\newcommand{\cG}{{\mathcal G}}
\newcommand{\cO}{{\mathcal O}}

\newcommand{\cN}{{\mathcal N}}

\newcommand {\QQ}{\mathbb{Q}}

\newcommand {\PP}{\mathbb{P}}

\newcommand {\YY}{\mathbb{Y}}
\newcommand {\VV}{\mathbb{V}}
\newcommand {\ra}{\longrightarrow}

\begin{document}
\title[]{Dimension of families of determinantal schemes.}

\author[Jan O.\ Kleppe, Rosa M.\ Mir\'o-Roig]{Jan O.\ kleppe, Rosa M.\
Mir\'o-Roig$^{*}$}
\address{Faculty of Engineering,
         Oslo University College,
         Cort Adelers gt. 30, N-0254 Oslo,
         Norway}
\email{JanOddvar.Kleppe@iu.hio.no}
\address{Facultat de Matem\`atiques,
Departament d'Algebra i Geometria, Gran Via de les Corts Catalanes
585, 08007 Barcelona, SPAIN } \email{miro@mat.ub.es}

\date{\today}
\thanks{$^*$ Partially supported by BFM2001-3584.}



\begin{abstract} A scheme $X\subset \PP^{n+c}$ of codimension $c$ is called
{\em standard determinantal} if its homogeneous saturated ideal
can be generated by the maximal minors of a homogeneous $t \times
(t+c-1)$ matrix and $X$ is said to be {\em good determinantal } if
it is standard determinantal and a generic complete intersection.
Given integers $a_0,a_1,...,a_{t+c-2}$ and $b_1,...,b_t$ we denote
by $W(\underline{b};\underline{a})\subset \Hi ^p(\PP^{n+c})$
(resp. $W_s(\underline{b};\underline{a})$) the locus of good
(resp. standard) determinantal schemes $X\subset \PP^{n+c}$ of
codimension $c$ defined by the maximal minors of a $t\times
(t+c-1)$ matrix $(f_{ij})^{i=1,...,t}_{j=0,...,t+c-2}$ where
$f_{ij}\in k[x_0,x_1,...,x_{n+c}]$ is a homogeneous polynomial of
degree $a_j-b_i$.

In this paper we address the following three fundamental problems
: To determine (1) the dimension of
$W(\underline{b};\underline{a})$ (resp.
$W_s(\underline{b};\underline{a})$) in terms of $a_j$ and $b_i$,
(2) whether the closure of $W(\underline{b};\underline{a})$ is an
irreducible component of $\Hi ^p(\PP^{n+c})$, and (3) when $\Hi
^p(\PP^{n+c})$ is generically smooth along
$W(\underline{b};\underline{a})$. Concerning question (1) we give
an upper bound for the dimension of
$W(\underline{b};\underline{a})$ (resp.
$W_s(\underline{b};\underline{a})$) which works for all integers
$a_0,a_1,...,a_{t+c-2}$ and $b_1,...,b_t$, and we conjecture that
this bound is sharp. The conjecture is proved for $2\le c\le 5$,
and for $c\ge 6$ under some restriction on $a_0,a_1,...,a_{t+c-2}$
and $b_1,...,b_t$. For questions (2) and (3) we have an
affirmative answer for $2\le c \le 4$ and $n\ge 2$, and  for $c\ge
5$ under certain numerical assumptions.
\end{abstract}


\maketitle

\tableofcontents


  \section{Introduction} \label{intro}

In this paper, we will deal with {\em determinantal} schemes, i.e.
schemes defined by the vanishing locus of the  minors of a
homogeneous polynomial matrix. Some classical schemes that can be
constructed in this way are the Segre varieties, the rational
normal scrolls and the Veronese varieties. Determinantal schemes
have been a central topic in both commutative algebra and
algebraic geometry and, due to their important role, their study
has attracted many researchers and have received considerable
attention in the literature.  Some of the most remarkable results
about determinantal schemes are due to J.A. Eagon and M. Hochster;
in \cite{e-h}, they proved that generic determinantal schemes are
arithmetically Cohen-Macaulay; and to J.A. Eagon and D.G.
Northcott; in \cite{e-n}, they constructed a finite free
resolution for any  standard determinantal and as a corollary they
got that standard determinantal schemes are arithmetically
Cohen-Macaulay. Since then many authors have made important
contributions to the study of determinantal schemes and the reader
can look at \cite{b-v}, \cite{dan}, \cite{BH} and \cite{eise} for background,
history  and a list of important papers.

\vskip 2mm A scheme $X\subset \PP^{n+c}$ of codimension $c$ is
called {\em standard determinantal} if its homogeneous saturated
ideal can be generated by the maximal minors of a homogeneous $t
\times (t+c-1)$ matrix and $X$ is said to be {\em good
determinantal} if it is standard determinantal and a generic
complete intersection. In this paper, we address the problem of
determining the dimension of the family of standard (resp. good)
determinantal schemes. The first important contribution to this
problem is due to G. Ellinsgrud \cite{elli}; in 1975, he proved
that every arithmetically Cohen-Macaulay, codimension 2 closed
subscheme $X$ of $\PP^{n+2}$ is unobstructed (i.e. the
corresponding point in the Hilbert scheme $\Hi ^p (\PP^{n+2})$ is
smooth) and he also computed the dimension of the Hilbert scheme
at $X$. Recall also that the homogeneous ideal of an
arithmetically Cohen-Macaulay, codimension 2 closed subscheme $X$
of $\PP^{n+2}$ is given by the maximal minors of a $t\times (t+1)$
homogeneous matrix, the Hilbert-Burch matrix. That is, such a
scheme is standard determinantal. The purpose of this work is to
extend Ellingsrud Theorem, viewed as a statement on standard
determinantal schemes of codimension 2, to arbitrary codimension.
The case of codimension 3, is solved in \cite{KMMNP}; Proposition
1.12. In the present  work, using essentially the methods
developed in \cite{KMMNP}; \S 10, we succeed in generalizing to
arbitrary codimension  the formula for the dimension of families
of determinantal schemes provided certain weak numerical
conditions are satisfied (see Theorem \ref{upperbound},
Proposition \ref{bound2} and Corollaries \ref{cod3}, \ref{cod4},
\ref{cod5}, \ref{cod6} and \ref{cod3dim0}). We also address the
problem whether the closure of the locus $W$ of determinantal
schemes in $\PP^{n+c}$ is an irreducible component of $\Hi
^p(\PP^{n+c})$ and when $\Hi ^p(\PP^{n+c})$ is generically smooth
along $W$ (see Corollaries \ref{component34}, \ref{component5},
\ref{component6} and \ref{c3c4}).

\vskip 2mm Next we outline the structure of the paper. In section
2, we recall the basic facts on standard and good determinantal
schemes $X\subset \PP^{n+c}$ of codimension $c$ defined by the
maximal minors of a $t\times (t+c-1)$ homogeneous matrix and the
associated complexes needed later on. Sections 3-5 are the heart
of the paper. Given integers $b_1,... ,b_t$ and $a_0, a_1 , ...
,a_{t+c-2}$,
 we  denote by
$W(\underline{b};\underline{a})\subset \Hi^{p}(\PP^{n+c})$  the
locus of good  determinantal schemes $X\subset \PP^{n+c}$ of
codimension $c\ge 2$ defined by the maximal minors of a
homogeneous matrix $\cA=(f_{ji})^{i=1,...,t}_{j=0,...,t+c-2}$
where $f_{ji}\in { k}[x_{0},...,x_{n+c}]$ is a homogeneous
polynomial of degree $a_j-b_{i}$. The goal of section 3 is to give
an upper bound  for the dimension of
$W(\underline{b};\underline{a})$ in terms of $b_1,...,b_t$ and
$a_0,a_1,...,a_{t+c-2}$ (cf. Theorem \ref{upperbound} and
Proposition \ref{bound2}). To this end we proceed by induction on
$c$ by successively deleting columns of the largest possible
degree and we strongly use the Eagon-Northcott complexes and the
Buchsbaum-Rim complexes associated to a standard determinantal
scheme. In section 4, using again induction on the codimension and
the theory of Hilbert flag schemes, we analyze when the upper
bound of $\dim W(\underline{b};\underline{a})$ given in section 3
is indeed the dimension of the determinantal locus. In turns out
that the upper bound of
 $\dim W(\underline{b};\underline{a})$ given in Theorem
 \ref{upperbound} is sharp in a number of instances. More
 precisely,
if $2\le c \le 3$, this is known (\cite{KMMNP}, \cite{elli}),  for
$4\le c \le 5$ it is a consequence of the main theorem of this
section (see Corollaries \ref{cod4} and \ref{cod5}), while for
$c\ge 6$ we  get the expected dimension formula for
$W(\underline{b};\underline{a})$ under more restrictive
assumptions (see Corollary \ref{cod6}).  In section 5, we study
when the closure of $ W(\underline{b};\underline{a})$ is an
irreducible component of $\Hi ^p(\PP^{n+c})$ and when $\Hi
^p(\PP^{n+c})$ is generically smooth along $
W(\underline{b};\underline{a})$, and other cases of
unobstructedness. The main result of this section (Theorem
\ref{main5}) shows that the closure of $
W(\underline{b};\underline{a})$ is a generically smooth
irreducible component provided the zero degree piece of certain
$\Ext^1$-groups vanishes. The  conditions of the Theorem can be
shown to be satisfied in a wide number of cases which we make
explicit through this section. In  particular we show that the
mentioned $\Ext^1$-groups vanishes if $3\le c \le 4$ (Corollary
\ref{component34}). Similarly, in Corollaries \ref{component5},
\ref{component6} and \ref{c3c4} and as a consequence of Theorem
\ref{main5}, we prove that under certain numerical assumptions the
closure of $ W(\underline{b};\underline{a})$ is indeed a
generically smooth, irreducible component of $\Hi ^p(\PP^{n+c})$
of the expected dimension. In Examples \ref{ex1} and \ref{ex2}, we
show that this is not always the case, although the examples
created are somewhat special because all the entries of the
associated matrix are linear entries and $
W(\underline{b};\underline{a})$ parameterizes curves.

We end the paper with a Conjecture raised by this paper and proved
in many cases (cf. Conjectures 6.1 and 6.2) and we correct an
inaccuracy in \cite{KMMNP}.

 \vskip 4mm The first author expresses  his thanks to the
 University of Barcelona and the University of Oslo.
 Part of this work was done while the second author
was a guest of the University of Oslo and she thanks the
University of Oslo for its hospitality.

 \vskip 4mm

{\bf Notation:} Throughout this paper $\PP^N$ will be the
$N$-dimensional projective space over an algebraically closed
field $k$, $R=k[x_0, x_1, \dots ,x_N]$ and $\goth m= (x_0, \dots
,x_N)$. The sheafification of a graded $R$-module $M$ will be
denoted by $\tilde{M}$.

For any closed subscheme $X$ of $\PP^N$ of codimension $c$, we
denote by  ${\mathcal I}_X$ its ideal sheaf,  ${\mathcal N}_X$ its
normal sheaf, $ I(X)=H^0_{*}({\mathcal I}_X)$ its saturated
homogeneous ideal and $\omega_X ={\mathcal E}xt^c_{{\mathcal
O}_{\PP^N}} ({\mathcal O}_X,{\mathcal O}_{\PP^N})(-N-1)$ its
canonical sheaf. For any quotient $A$ of $R$ of codimension $c$,
we let $I_A=\ker(R\twoheadrightarrow A)$, $N_A=\Hom_R(I_A,A)$ be
the normal module and $K_A=\Ext^c_R (A,R)(-N-1)$ be its canonical
module. When we write $X=Proj(A)$, we let $A=R/I(X)$ and
$K_X=K_A$.

We denote the Hilbert scheme by $\Hi ^p(\PP^N)$. Thus, any $X\in
\Hi ^p(\PP^N)$ is a closed subscheme of $\PP^N$ with Hilbert
polynomial $p\in \QQ[s]$. By definition $X\in \Hi ^p(\PP^N)$ is
unobstructed if $\Hi ^p(\PP^N)$ is smooth at $X$.


\section{Preliminaries}

This section provides the background and  basic results on
standard (resp. good) determinantal schemes needed in the sequel
and we refer to \cite{b-v} and \cite{eise} for more details.

\vskip 2mm

Let $\cA$ be a homogeneous matrix, i.e. a matrix representing a
degree 0 morphism $\phi $ of free graded $R$-modules. In this
case, we denote by $I(\cA)$ (or $I(\phi )$)  the ideal of $R$
generated by the maximal minors of $\cA$.

\begin{definition}
 A codimension $c$ subscheme $X\subset \PP^{n+c}$ is
called a \emph{standard determinantal} scheme if $I(X)=I(\cA)$ for
some $t\times (t+c-1)$ homogeneous matrix $\cA$.  $X\subset
\PP^{n+c}$ is called a \emph{good determinantal} scheme if
additionally, $\cA$ contains a $(t-1)\times (t+c-1)$ submatrix
(allowing a change of basis if necessary) whose ideal of maximal
minors define a scheme of codimension $c+1$.
 \end{definition}

It is well known that a good determinantal scheme $X\subset
\PP^{n+c}$ is standard determinantal and the converse is true
provided $X$ is a generic complete intersection, cf. \cite{KMNP}.

\vskip 2mm
 Now we are going to describe
the generalized Koszul complexes associated to a codimension $c$
standard determinantal scheme $X$. To this end, we denote by
$\varphi :F \longrightarrow G$ the morphism of free graded
$R$-modules of rank $t$ and $t+c-1$, defined by the homogeneous
matrix $\cA$ of $X$. We denote by ${\cC}_i(\varphi ^*)$ the
generalized Koszul complex:

$${\cC}_i(\varphi^*): \  0 \longrightarrow \wedge^{i}G ^*\otimes
S_{0}(F^*)\ra \wedge^{i-1} G^* \otimes S _{1}( F^*)\ra \ldots \ra
\wedge ^{0} G^* \otimes S_i (F^*) \ra 0.$$

 Let ${\cC}_i(\varphi^*)^*$ be the $R$-dual of
${\cC}_i(\varphi^*)$. The  map $\varphi$ induces graded morphisms
$$\mu_i:\wedge ^{t+i}G^*\otimes \wedge^tF\ra\wedge^{i}G^*.$$

They can be used to splice the complexes ${\cC}_{c-i-1}(\varphi ^*
)^*\otimes \wedge^{t+c-1}G^*\otimes \wedge^tF$ and
 ${\cC}_i(\varphi^*)$ to a complex ${\cD}_i(\varphi^*):$

$$0 \ra \wedge^{t+c-1}G^* \otimes S_{c-i-1}(F)\otimes \wedge^tF\ra
\wedge^{t+c-2} G ^*\otimes S _{c-i-2}(F)\otimes \wedge ^tF\ra
\ldots \ra$$ $$\wedge^{t+i}G^* \otimes S_{0}(F)\otimes
\wedge^tF\ra \wedge^{i} G ^*\otimes S _{0}(F^*) \ra \wedge ^{i-1}
G^* \otimes S_1(F^*)\ra \ldots \ra \wedge^0G^*\otimes S_i(F^*)\ra
0 .$$

\vskip 4mm The complex  ${\cD}_0(\varphi^*)$ is called the {\em
Eagon-Northcott complex} and the complex  ${\cD}_1(\varphi^*)$ is
called the {\em Buchsbaum-Rim complex}. Let us rename the complex
${\cC}_c(\varphi^*)$ as  ${\cD}_c(\varphi^*)$. Then we have the
following well known result:

\begin{proposition} \label{resolut} Let $X\subset \PP^{n+c}$ be a standard determinantal subscheme
of codimension $c$ associated to a graded minimal (i.e.
$\im(\varphi )\subset \mathfrak{m}G$) morphism $\varphi: F\ra G$
of free $R$-modules of rank $t$ and $t+c-1$, respectively. Set $M=
\coker (\varphi^*)$. Then it holds:

(i) ${\cD}_i(\varphi^*)$ is acyclic for $-1\le i \le c$.

(ii) ${\cD}_0(\varphi^*)$ is a minimal free graded $R$-resolution
of $R/I(X)$  and  ${\cD}_i(\varphi^*)$ is a minimal free graded
$R$-resolution of length $c$ of $S_i(M)$, $1\le i \le c$.

 (iii) $K_X \cong
S_{c-1}(M)$ up to degree shift . So, up to degree shift,
$\cD_{c-1}(\varphi^*)$ is a minimal free graded $R$-module
resolution of $K_X$.

\end{proposition}

\begin{proof}
See, for instance \cite{b-v}; Theorem 2.20  and \cite{eise};
Corollary A2.12 and Corollary A2.13.
\end{proof}

\vskip 2mm
\begin{remark} {\rm By Proposition \ref{resolut}(ii), any standard determinantal
scheme $X\subset \PP^{n+c}$ is arithmetically Cohen-Macaulay
(briefly, ACM).
 Moreover, in codimension 2, the converse
is true: If $X\subset \PP^{n+2}$ is an ACM, closed subscheme of
codimension 2 then it is standard determinantal (Hilbert-Burch
Theorem).}
\end{remark}

The homogeneous matrix $\cA$ associated to a standard
determinantal scheme $X\subset \PP^{n+c}$ of codimension $c$ also
defines an injective morphism $\varphi :\cF \longrightarrow \cG $
of locally free ${\cO}_{\PP^{n+c}}$-modules of rank $t$ and
$t+c-1$. Since the construction of the generalized Koszul
complexes globalizes, we can also associated to $\varphi ^*$ the
{\em Eagon-Northcott complex} of $\cO _{\PP^{n+c}}$-modules

$$0 \ra \wedge^{t+c-1}\cG^* \otimes S_{c-1}(\cF )\otimes
\wedge^t\cF\ra \wedge^{t+c-2} \cG^* \otimes S _{c-2}(\cF )\otimes
\wedge ^t\cF\ra \ldots \ra$$ $$\wedge^{t}\cG^* \otimes
S_{0}(\cF)\otimes \wedge^t\cF \ra {\cO}_{\PP^{n+c}} \ra {\cO}_X\ra
0 $$

\noindent and the {\em Buchsbaum-Rim complex} of locally free $\cO
_{\PP^{n+c}}$-modules

$$0 \ra \wedge^{t+c-1}\cG^* \otimes S_{c-2}(\cF )\otimes
\wedge^t\cF \ra \wedge^{t+c-2} \cG^* \otimes S _{c-3}(\cF )\otimes
\wedge ^t\cF \ra \ldots \ra$$ $$\wedge^{t+1}\cG^* \otimes
S_{0}(\cF)\otimes \wedge^t\cF \ra \cG^* \stackrel {\varphi ^* }{
\longrightarrow}  \cF^* \ra \tilde{M} \ra 0. $$

Since the degeneracy locus of $\varphi ^*$ has codimension $c$,
these two complexes are acyclic. Moreover, the kernel of $\varphi
^*$ is called the {\em 1st Buchsbaum-Rim sheaf} associated to
$\varphi ^*$.

\vskip 2mm
 Let $X\subset \PP^{n+c}$ be a standard (resp. good) determinantal scheme
 of codimension    $c\ge 2$
defined by the vanishing of the maximal minors of a $t\times
(t+c-1)$ matrix $\cA=(f_{ji})_{i=1,...t}^{j=0,...,t+c-2}$ where
$f_{ji}\in { k}[x_{0},...,x_{n+c}]$ are homogeneous polynomials of
degree $a_j-b_{i}$ with $b_1 \le ... \le b_t$ and $a_0 \le a_1\le
... \le a_{t+c-2}$. We assume without loss of generality that
$\cA$ is minimal; i.e., $f_{ji}=0$ for all $i,j$ with
$b_{i}=a_{j}$. If we let $u_{ji}=a_j-b_i$ for all $j=0, \dots ,
t+c-1$ and $i=1, \dots . t$, the matrix
$\cU=(u_{ji})_{i=1,...t}^{j=0,...,t+c-2}$ is called the {\em
degree matrix} associated to $X$. We have:

\begin{lemma} \label{columns} The matrix $\cU$ has the following properties:
\begin{itemize}
\item[(i)] For every $j$ and $i$, $u_{j,i}\le u_{j+1,i}$ and $u_{j,i}\ge
u_{j,i+1}$.
\item[(ii)]
For every $i=1,\dots , t$, $u_{i-1,i}=a_{i-1}-b_i> 0$.
\end{itemize}
And vice versa, given a degree matrix $\cU $ of integers verifying
(i) and (ii) there exists a codimension $c$ standard (resp. good)
determinantal scheme $X\subset \PP^{n+c}$ with associated degree
matrix $\cU .$
\end{lemma}
\begin{proof} The first condition is obvious. For the second one
we only need to observe that if for some $i=1, \dots , t$, we have
$u_{i-1,i}\le 0$ then in the matrix $\cA$ we have $f_{j,k}=0$ for
$j\le i-1$ and $k\ge i$. But this would imply that the minor which
is obtained by deleting the last $c-1$ columns has to be zero
contradicting the minimality of $\cA$.

The converse is trivial. Indeed, given a matrix of integers, $\cU
$,
 satisfying (i) and (ii), we can consider the standard (resp.
 good)
determinantal scheme $X\subset \PP^{n+c}$ of codimension $c$
associated to the homogeneous matrix

$$ \cA=\begin{pmatrix} x_{0}^{a_{t+c-2}-b_t} & x_{1}^{a_{t+c-3}-b_t} &
... & ... & x_{c-1}^{a_{t-1}-b_{t}}  & 0 & 0 & ...& ...
\\ 0 & x_{0}^{a_{t+c-3}-b_{t-1}} & x_{1}^{a_{t+c-4}-b_{t-1}} & ... & ... &
x_{c-1}^{a_{t-2}-b_{t-1}}  & 0  & 0 & ... \\ 0 & 0 &
x_{0}^{a_{t+c-4}-b_{t-2}} & x_{1}^{a_{t+c-5}-b_{t-2}} & ... & ...
& x_{c-1}^{a_{t-3}-b_{t-2}}  & 0  & ...
\\ .... & ... & ... & ... & ... & ...  & ... & ... & ...
\end{pmatrix}
$$

(resp. $ \cA=$ $$\begin{pmatrix} x_{0}^{a_{t+c-2}-b_t} &
x_{1}^{a_{t+c-3}-b_t} & ... & ... & x_{c-1}^{a_{t-1}-b_{t}}  & 0 &
0 & ...& ...
\\ x_c^{a_{t+c-2}-b_{t-1}} & x_{0}^{a_{t+c-3}-b_{t-1}} & x_{1}^{a_{t+c-4}-b_{t-1}} & ... & ... &
x_{c-1}^{a_{t-2}-b_{t-1}}  & 0  & 0 & ... \\ 0 &
x_c^{a_{t+c-3}-b_{t-2}} & x_{0}^{a_2-b_{t-2}} &
x_{1}^{a_{t+c-5}-b_{t-2}} & ... & ... & x_{c-1}^{a_{t-3}-b_{t-2}}
& 0 & ...
\\ .... & ... & ... & ... & ... & ...  & ... & ... & ...
\end{pmatrix} ).
$$ Up to re-ordering, we easily check that the degree matrix
associated to $X$ is $\cU$.
\end{proof}

Given integers $b_1,...,b_{t}$ and $a_0, a_1, ... , a_{t+c-2}$, we
denote by $W(\underline{b};\underline{a})\subset
\Hi^{p}(\PP^{n+c})$ (resp.
$W_s(\underline{b};\underline{a})\subset \Hi^{p}(\PP^{n+c})$) the
locus of good (resp. standard) determinantal schemes $X\subset
\PP^{n+c}$ of codimension $c\ge 2$ defined by the maximal minors
of a homogeneous matrix $\cA=(f_{ji})_{j=0,...,t+c-2}^{i=1,...,t}$
where $f_{ji}\in { k}[x_{0},...,x_{n+c}]$ is a homogeneous
polynomial of degree $a_j-b_{i}$. Clearly,
$W(\underline{b};\underline{a})\subset
W_s(\underline{b};\underline{a}).$ Moreover, it holds:

\begin{corollary} \label{WWs} Assume $b_1\le ... \le b_t$
and $a_0\le a_1 \le ... \le a_{t+c-2}$. We have that
$W(\underline{b};\underline{a})\ne \emptyset $ if and only  if
$W_s(\underline{b};\underline{a}) \ne \emptyset $ if and only if
$u_{i-1,i}=a_{i-1}-b_i> 0$ for $i=1,...,t$.
\end{corollary}
\begin{proof} It easily follows from Lemma \ref{columns}.
\end{proof}

Let $X\subset \PP^{n+c}$ be a  good scheme
 of codimension    $c\ge 2$
defined by the homogeneous  matrix
$\cA=(f_{ji})_{i=1,...t}^{j=0,...,t+c-2}$. It is well known that
by successively deleting columns from the right hand side, one
gets a flag of determinantal subschemes
\begin{equation}\label{flag}
 ({\mathbf X.}) :  X = X_c \subset X_{c-1} \subset ...  \subset X_{2} \subset
\PP^{n+c}
\end{equation}
where each $X_{i+1} \subset X_i$ (with ideal sheaf ${\mathcal
I_{X_{i+1}|X_i}} = {\mathcal I_i}$) is of codimension 1, $X_{i}
\subset \PP^{n+c}$ is of codimension $i$ ($i=2,\dots , c$) and
where there exist ${\mathcal O}_{X_i}$-modules ${\mathcal M_i}$
fitting into short exact sequences
\begin{equation}\label{fix-ref}
0\rightarrow {\mathcal O}_{X_i}(-a_{t+i-1})\rightarrow {\mathcal
M}_i \rightarrow {\mathcal M}_{i+1} \rightarrow 0 \ \  {\rm for} \
\ 2 \leq i \leq c-1,
\end{equation}
 such that ${\mathcal I}_i (a_{t+i-1})$
is the ${\mathcal O}_{X_i}$-dual of ${\mathcal M}_i$,  for $2 \leq
i \leq c$, and ${\mathcal M}_{2}$ is a twist of the canonical
module of $X_2$, cf. (3.4)-(3.7) for details.

\begin{remark} \label{dep} Assume $b_1\le ... \le b_t$ and
$a_0\le a_1 \le ... \le a_{t+c-2}$. Similar  to Lemma
\ref{columns}, if $X$ is general in
$W(\underline{b};\underline{a})$ and $u_{i-min(\alpha
,t),i}=a_{i-min(\alpha ,t)}-b_i\ge 0$ for $min(\alpha ,t)\le i\le
t$, then $X_{j}=Proj(D_j)$, for all $j=2, \cdots , c$, is
non-singular except for a subset of codimension $\ge min\{2\alpha
-1, j+2 \}$, i.e.
\begin{equation}\label{a-b} \codim_{X_j}Sing(X_j)\ge min\{2\alpha
-1, j+2\} . \end{equation} This follows from \cite{chang}; Theorem
arguing as in \cite{chang}; Example 2.1. In particular, if $\alpha
\ge 3$, we get that for each $i>0$, the closed embeddings
$X_i\subset \PP^{n+c}$ and $X_{i+1}\subset X_{i}$ are local
complete intersections outside some set $Z_i$ of codimension $\ge
min(4,i+1)$ in $X_{i+1}$ ($\depth_{Z_i}\cO_{X_{i+1}}\ge
\mbox{min}(4,i+1)$).
\end{remark}


\section{Upper bound for the dimension of the determinantal locus}
 \label{upper bound}

The goal of this section is to write down an upper bound for the
dimension of the locus $W(\underline{b};\underline{a})$ (resp.
$W_s(\underline{b};\underline{a})$) of good (resp. standard)
determinantal subschemes $X \subset \proj{n+c}$  of codimension
$c$ inside the Hilbert scheme $\Hi^{p(s)} (\PP^{n+c} )$, where
$p(s)\in \QQ [s]$ is the Hilbert polynomial of $X$. In  section 4,
we will analyze when the mentioned upper bound is sharp and in
section 5, we will discuss under which conditions the closure of
$W$ in $\Hi^{p(s)} (\PP^{n+c} )$ is a generically smooth,
irreducible component of $\Hi ^{p(s)} (\PP^{n+c} )$.

\vskip 4mm Let $X\subset \PP^{n+c}$ be a good  determinantal
scheme
 of codimension    $c\ge 2$
defined by the vanishing of the maximal minors of a $t\times
(t+c-1)$ matrix $\cA=(f_{ji})_{i=1,...t}^{j=0,...,t+c-2}$ where
$f_{ji}\in { k}[x_{0},...,x_{n+c}]$ are homogeneous polynomials of
degree $a_j-b_{i}$ and let $A=R/I(X)$ be the homogeneous
coordinate ring of $X$. The matrix $\cA$ defines a morphism of
locally free sheaves
\[
\varphi :{\mathcal F}:=\bigoplus_{i=1}^{t}{\mathcal
O}_{\proj{n+c}} (b_i)\longrightarrow {\mathcal
G}:=\bigoplus_{j=0}^{t+c-2}{\mathcal O}_{\proj{n+c}} (a_j)
\]
and we may assume without loss of generality that $\varphi $ is
minimal; i.e., $f_{ji}=0$ for all $i,j$ with $b_{i}=a_{j}$.

 \vskip 2mm
Our aim is to upper bound $\dim W(\underline{b};\underline{a})$ in
terms of $b_1,...,b_t$ and $a_0,a_1,...,a_{t+c-2}$. To this end,
we consider the affine scheme $\VV =\Hom_{{\mathcal
O}_{\PP^{n+c}}}({\mathcal F},{\mathcal G})$
 whose rational points are the
morphisms from ${\mathcal F}$ to ${\mathcal G}$.
 Let $\YY $ be the
non-empty, open, irreducible subscheme of $\VV $ whose rational
points are the morphisms $\varphi _{\lambda }:{\mathcal
F}\longrightarrow {\mathcal G}$ such that their associated
homogeneous matrix $\cA _{\lambda }$ defines a good determinantal
subscheme $X_{\lambda }\subset \PP^{n+c}$.

\vskip 4mm The Eagon-Northcott complex of the universal morphism
$$\Psi: pr_{2}^{*}{\mathcal F}\longrightarrow pr^{*}_{2}{\mathcal
G}$$ \noindent on $\YY \times \PP ^{n+c}$ (where $pr_{2}:\YY\times
\PP^{n+c}\longrightarrow \PP^{n+c}$ is the natural projection)
induces a morphism

$$f:\YY\longrightarrow W(\underline{b};\underline{a})$$

\noindent which is defined by $f(\varphi _{\lambda }):=X_{\lambda
}$ on closed points. The affine group scheme $G:=Aut({\mathcal
F})\times Aut({\mathcal G})$ operates on $\YY$:

$$ \sigma :G\times \YY\longrightarrow \YY; \quad ((\alpha ,\beta
),\varphi _{\lambda }) \mapsto \beta \varphi _{\lambda
}\alpha^{-1}. $$

The action $\sigma $ is compatible with the morphism $f$. Thus, at
least set-theoretically $f:\YY\longrightarrow
W(\underline{b};\underline{a})$ induces a surjective map from the
orbit set $\YY //G$ to $W(\underline{b};\underline{a})$. Moreover,
since the map from $\YY$ to the closure
$\overline{W(\underline{b};\underline{a})}$ in $\Hi
^{p(s)}(\proj{n+c})$ is dominant, we get that
$W(\underline{b};\underline{a})$ is irreducible and
 we have (small letters denote dimension):
\begin{equation} \label{rosa-(1)}
\dim W(\underline{b};\underline{a}) \le \hom_{{\mathcal
O}_{\PP^{n+c}}}({\mathcal F},{\mathcal G})- \aut({\mathcal
G})-\aut({\mathcal F})+ \dim (G_{\lambda })
\end{equation}
where
\[
G_{\lambda }  =  \{(\delta ,\tau )\in Aut({\mathcal F})\times
Aut({\mathcal G}) \mid  \tau \varphi _{\lambda }\delta
^{-1}=\varphi _{\lambda }\}
\]
 is the isotropy group of a
general closed point $\varphi _{\lambda }\in \YY$. By
\cite{KMMNP}; Proposition 10.2, for all $\varphi _{\lambda }\in
\YY $,
 we have (we let $\binom{n+a}{ a }=0$ for $a<0$, as usual):
\begin{equation}
\dim (G_{\lambda }) =  \aut( {\mathcal B}_{\lambda }) +\sum_{j,i}
\binom{b_i-a_j+n+c}{n+c}.
\end{equation}

Therefore, we have
\begin{equation}\label{boundw}
 \dim W(\underline{b};\underline{a}) \le \sum_{i,j}
\binom{a_i-b_j+n+c}{n+c}- \sum _{i,j} \binom{a_i-a_j+n+c}{n+c}-
 \end{equation}
\[
 \sum _{i,j} \binom{b_i-b_j+n+c}{n+c} + \sum _{j,i}
\binom{b_j-a_i+n+c}{n+c}+\aut({\mathcal B}_{\lambda }).
\]

 Our next goal is to bound $ \aut( {\mathcal
B})$ in terms of $a_j$ and $b_i$, where ${\mathcal
B}=\coker(\varphi )$ and  $\varphi $ is a closed point of $\YY $.
To this end we need to fix some more notation.

\vskip 2mm Let ${\cA}_i$ be the matrix obtained deleting the last
$c-i$ columns. The matrix ${\cA}_i$ defines a morphism

\begin{equation}\label{gradedmorfismo} \varphi
_i:F=\bigoplus _{i=1}^tR(b_i)\longrightarrow G_i:=\bigoplus
_{j=0}^{t+i-2}R(a_j)
\end{equation}
of $R$-free modules and let $B_i$ be the cokernel of $\varphi _i$.
 Put $\varphi =\varphi _0 $, $G=G_c$ and $B=B_c$. Let
$M_i$ be the cokernel of $\varphi _i^*=\Hom_R(\varphi _i,R)$, i.e.
let the  sequence

\begin{equation}\label{defMi}
G_i^*\stackrel {\varphi_i^*}{ \longrightarrow} F^* \longrightarrow
M_i\cong  \Ext ^1_R(B_i,R)\longrightarrow 0
\end{equation}
be exact. If $D_i\cong R/I_{D_i}$ is the $k$-algebra given by the
maximal minors of ${\cA}_i$ and $X_i=Proj(D_i)$ (i.e.
$R\twoheadrightarrow D_2 \twoheadrightarrow D_3 \twoheadrightarrow
...\twoheadrightarrow D_c=A$), then $M_i$ is a $D_i$-module and
there is an exact sequence

\begin{equation}\label{Mi}
0\longrightarrow D_i \longrightarrow M_i(a_{t+i-1})
\longrightarrow M_{i+1}(a_{t+i-1}) \longrightarrow 0
\end{equation}
in which $D_i \longrightarrow M_i(a_{t+i-1})$ is the regular
section which defines $D_{i+1}$ \cite{KMNP}. Indeed,

\begin{equation}\label{Di}
0\longrightarrow
M_i(a_{t+i-1})^*=\Hom_{D_i}(M_i(a_{t+i-1}),D_i)\longrightarrow D_i
\longrightarrow D_{i+1}\longrightarrow 0
\end{equation}
and we may put $I_i:=I_{D_{i+1}/D_i}=M_i(a_{t+i-1})^*$. An
$R$-free resolution of $M_i$ is given by Proposition 2.2, and we
get in particular that $M_i$ is a maximal Cohen-Macaulay
$D_i$-module. Using (\ref{Di}) we see that $I_i$ is also a maximal
Cohen-Macaulay $D_i$-module. Proposition 2.2 (iii) also gives us
$K_{D_i}(n+c+1)\cong S_{i-1}M_i(\ell _i)$ where $\ell
_i:=\sum_{j=0}^{t+i-2}a_j-\sum_{i=1}^tb_i$.

\vskip 2mm In what follows we always let $Z_i\subset X_i$ be some
closed subset such that $U_{i}=X_{i}-Z_{i}\hookrightarrow
\PP^{n+c}$ is a local complete intersection. By the well known
fact that the 1. Fitting ideal of $M_{i}$ is equal to
$I_{t-1}(\varphi _{i})$, we get that $\tilde{M}_{i}$ is locally
free of rank 1 precisely on $X_i-V(I_{t-1}(\varphi _{i}))$
\cite{BH}, Lemma 1.4.8. Since the set of non locally complete
intersection points of $X_i\hookrightarrow \PP^{n+c}$ is precisely
$V(I_{t-1}(\varphi _i))$ by e.g. \cite{ulr}, Lemma 1.8, we get
that $U_i\subset X_i-V(I_{t-1}(\varphi_{i}))$ and that
$\tilde{M}_{i}$ and  $\cI _{X_{i}}/\cI ^2_{X_{i}}$  are locally
free on $U_i$.

Finally note that there is a close relation between $
M_{i+1}(a_{t+i-1})$ and the normal module
$N_{D_{i+1}/D_i}:=\Hom_{D_i}(I_{i},D_{i+1})$ of the quotient
$D_i\rightarrow D_{i+1}$. If we suppose $\depth_{I(Z_i)}D_i\ge 2$,
we get, by applying $ \Hom_{D_i}(I_i,.)$ to (\ref{Di}), that

\begin{equation} \label{DiMi}
0\longrightarrow D_i \longrightarrow M_i(a_{t+i-1})
\longrightarrow N_{D_{i+1}/D_i}
\end{equation}
is exact. Hence we have an injection
$M_{i+1}(a_{t+i-1})\hookrightarrow N_{D_{i+1}/D_i}$, which in the
case $\depth_{I(Z_i)}D_i\ge 3$ leads to an isomorphism
$M_{i+1}(a_{t+i-1})\cong N_{D_{i+1}/D_i}$. Indeed, this follows
from the more general fact  (by letting $M=N=I_i$) that if $M$ and
$N$ are finitely generated $D$-modules such that
$\depth_{I(Z)}M\ge r+1$ and $\tilde{N}$ is locally free on
$U:=X-Z$ ($X=Proj(D)$), then the natural map

\begin{equation} \label{NM}
\Ext^{i}_D(N,M)\longrightarrow
H_{*}^{i}(U,{\cH}om_{{\cO}_X}(\tilde{N},\tilde{M}))\simeq
H^{i+1}_{I(Z)}(\Hom_D(N,M))
\end{equation}
is an isomorphism, (resp. an injection) for $i<r$ (resp. $i=r$),
Cf. \cite{SGA2}, exp. VI.

\begin{lemma} \label{key}
Let $M$ be an $R$-module. With the above notation, the sequence
$$0\rightarrow \Hom_R(M_i,M)\rightarrow F\otimes_R M \rightarrow
G_i\otimes _R M\rightarrow B_i\otimes _R M\rightarrow 0$$ is exact
and $\Hom_R(M_i,M)=\TT_1^R(B_i,M)$.

\end{lemma}

\begin{proof} We apply $\Hom(.,R)$ to
\begin{equation}\label{Bi}
0 \longrightarrow F\longrightarrow G_i\longrightarrow B_i
\longrightarrow 0 \end{equation} and we get $$ 0\rightarrow
\Hom(B_i,R)\rightarrow G_i^*\rightarrow F^*\rightarrow
\Ext^1_R(B_i,R)=M_i \rightarrow 0.$$ Hence $$0\rightarrow
\Hom(M_i,M)\rightarrow \Hom(F^*,M)\cong F\otimes M\rightarrow
\Hom(G_i^*,M)\cong G_i\otimes M $$ and we get the first exact
sequence and $\Hom(M_i,M)=\TT_1^R(B_i,M)$ by applying $(.)\otimes
_R M$ to (\ref{Bi}).
\end{proof}

\begin{lemma}\label{key2} With the notations above, if $C$ is good
determinantal, then  $\depth _{I(Z_i)}D_i\ge 1$ for $2\le i \le c$
and $\Hom _{D_{i}}(M_i,M_i)=D_i.$
\end{lemma}
\begin{proof}
If $X$ is a standard determinantal scheme, defined by some
$t\times (t+i-1)$ matrix, and if we delete a column and let $Y$ be
the corresponding determinantal scheme, then $Y$ is also standard
determinantal \cite{B}. Hence if $X$ is good determinantal, it
follows that $Y$ is also good determinantal by the definition of a
good determinantal scheme. In particular all $X_i$, $2\le i \le
c$, are good determinantal schemes and hence generic complete
intersections. By the definition of $Z_i$, we get $\depth
_{I(Z_{i})}D_i\ge 1$.

 Let $U_i=Proj(D_i)-Z_i$ and note that $\tilde{M_i}|_{U_i}$ is
an invertible sheaf. Let $S_r(M_i)$ be the $r$-th symmetric power
of the $D_i$-module $M_i$. For $1\le r \le i-1$, $S_r(M_i)$ are
maximal Cohen-Macaulay modules and $S_{i-1}(M_i)(\ell _i)\cong
K_{D_i}(n+c+1)$ (cf. Proposition 2.2 (iii)). By (\ref{NM}) we have
injections $$\Hom_{D_i}(S_rM_i,S_rM_i)\hookrightarrow
H^0_{*}(U_i,{\cH}om(\tilde{S_rM_i},\tilde{S_rM_i}))\cong
H^0_{*}(U_i,\tilde{D_i}).$$

Since $S_{i-1}(M_i)$ is a twist of the canonical module  and since
$\Hom (K_{D_i},K_{D_i})\cong D_i$, we get the lemma from the
commutative diagram

\[
\begin{array}{cccccc}

\Hom_{D_i}(M_i,M_i) &  \hookrightarrow & H^0_{*}(U_i,\tilde{D_i})
\\ \psi \downarrow & &
\parallel  &     \\ \Hom_{D_i}(S_{i-1}M_i,S_{i-1}M_i) &
 \hookrightarrow & H^0_{*}(U_i,\tilde{D_i}).
\end{array}
\]
Indeed, $\psi $ is injective and we conclude by $$D_i\rightarrow
\Hom_{D_i}(M_i,M_i)\hookrightarrow
\Hom_{D_i}(S_{i-1}M_i,S_{i-1}M_i)\cong D_i.$$
\end{proof}

\begin{proposition} \label{aut(B)} Assume $b_1\le ... \le b_t$,
$a_0\le a_1 \le ... \le a_{t+c-2}$ and $c\ge 2$.  Set $\ell
:=\sum_{j=0}^{t+c-2}a_j-\sum_{i=1}^tb_i$. If $(c-1)a_{t+c-2}<\ell
$ then, $\aut(\cB )=1$. Otherwise, it holds $$\aut(\cB)\le
\sum_{i=1}^{c-3}\left (\sum _{r+s=i} \sum _{ 0\le i_i<...<i_r\le
t+i \atop 1\le j_1\le ... \le j_s \le t} (-1)^{i-r}
\binom{h_i+a_{i_1}+\cdots +a_{i_r}+b_{j_1}\cdots +b_{j_s} }{n+c}
\right )$$

$$+\binom{h_0}{n+c} +1 $$ \noindent where we set $h_i:=
2a_{t+1+i}+a_{t+2+i}+\cdots a_{c+t-3}+a_{c+t-2}-\ell +n+c$, for
$i=0,1,\cdots , c-3$.
\end{proposition}

\begin{proof}  If $(c-1)a_{t+c-2}< {\ell }$ then, by
\cite{KMMNP}; Lemma 10.1(ii), ${\mathcal B}$ is stable and $ \aut(
{\mathcal B})=1$ because stable reflexive sheaves are simple. So,
from now on, we assume $(c-1)a_{t+c-2}\ge  {\ell }$ and we will
proceed by induction on $c$ by successively deleting columns from
the right side, i.e.  of the largest degree. For $c=2$ the result
was proved in \cite{elli} if $n\ge 1$ and in \cite{KMMNP} for any
$n\ge 0$. So, we will assume $c\ge 3$.

Consider the commutative diagram
\[
\begin{array}{cccccccccccc}
 & &  &  & 0 &  & 0 & & &
\\
&& && \downarrow  &&   \downarrow  \\ & &  &  & \cO
_{\PP^{n+c}}(a_{t+c-2}) & = & \cO _{\PP^{n+c}}(a_{t+c-2}) & & &
\\ && &&
\downarrow  &&   \downarrow  \\ 0 & \longrightarrow &
\oplus_{i=1}^t\cO _{\PP^{n+c}}(b_i) & \stackrel {\varphi_c}{
\longrightarrow} & \oplus_{j=0}^{t+c-2}\cO _{\PP^{n+c}}(a_j) &
\longrightarrow & {\mathcal B}_{c} & \longrightarrow & 0
\\ && \| &&
\downarrow  &&   \downarrow  \\
 0 & \longrightarrow & \oplus_{i=1}^t\cO _{\PP^{n+c}}(b_i) &
\stackrel {\varphi_{c-1}}{ \longrightarrow} &
\oplus_{j=0}^{t+c-3}\cO _{\PP^{n+c} }(a_j) & \longrightarrow &
{\mathcal B}_{c-1 } & \longrightarrow & 0 \\
 && && \downarrow  &&   \downarrow  \\& &  &  & 0 &  & 0 & & &
\end{array}
\]
and the exact sequence $$ 0 \ra \Hom({\mathcal B}_{c-1 },
{\mathcal B}_{c}) \ra Aut( {\mathcal B}_{c })\stackrel {\alpha }
{\longrightarrow} \Hom(\cO_{\PP^{n+c}}(a_{t+c-2}), {\mathcal
B}_{c})=H^0 ({\mathcal B}_{c}(-a_{t+c-2})) \ra $$ $$ \Ext^1_{\cO
_{\PP^{n+c}}}(\cB_{c-1},\cB_c)  \ra \Ext^1_{\cO
_{\PP^{n+c}}}(\cB_{c},\cB_c) \ra 0.$$

Moreover, if we tensor with $.\otimes_R B_c$ the exact sequence
$$0\rightarrow D_ {c-1} (-a_{t+c-2})\ra  M_{c-1}\ra M_c\ra 0$$ we
get $$ \TT^R_1(B_c,M_{c-1})\ra \TT_1^R(B_c,M_c) \ra
D_{c-1}(-a_{t+c-2})\otimes  B_c\ra $$ $$\ra  M_{c-1}\otimes
B_c\cong \Ext^1_R(B_{c-1},B_c) \ra M_c\otimes B_c\cong
\Ext^1_R(B_{c},B_c) \ra 0.$$

\vskip 2mm Applying Lemmas \ref{key} and \ref{key2} we get
$\TT^R_1(B_c,M_{c-1})=\Hom(M_c,M_{c-1})=0$ (since $M_c$ is
supported in  $X_c$ which has codimension 1 in
$X_{c-1}=Supp(M_{c-1})$ ) and
$\TT_1^R(B_c,M_c)=\Hom(M_c,M_c)=D_c=A.$ Hence,
$$H^0(\cB_c(-a_{t+c-2})) \rightarrow \Ext^1_{\cO
_{\PP^{n+c}}}(\cB_{c-1}, \cB_c)$$ coincides with
$(D_{c-1}(-a_{t+c-2})\otimes B_c)_0\rightarrow (M_{c-1}\otimes
B_c)_0$ whose kernel is $A_0\cong k$, i.e. 1-dimensional.
Therefore, $\dim (\im(\alpha))=1$ which gives us
\begin{equation}\label{aut1} \aut( {\mathcal
B}_{c})\le \hom({\mathcal B}_{c-1 }, {\mathcal B}_{c })+1.
\end{equation}

We call $0\ne e \in \Ext^1({\mathcal B}_{c-1
},\cO_{\PP^{n+c}}(a_{t+c-2}))$ the non-trivial extension $$ e:
\quad 0 \ra \cO_{\PP^{n+c}}(a_{t+c-2}) \ra {\mathcal B}_{c } \ra
{\mathcal B}_{c-1} \ra 0, $$ we have $$ 0 \ra \Hom({\mathcal
B}_{c-1 }, \cO(a_{t+c-2})) \ra  \Hom({\mathcal B}_{c-1 },
 {\mathcal B}_{c})\stackrel {\eta}{ \longrightarrow}
 \Hom({\mathcal B}_{c-1
},{\mathcal B}_{c-1 } ) \stackrel {\delta}{ \longrightarrow}
\Ext^1({\mathcal B}_{c-1 }, \cO(a_{t+c-2})) .$$ Since $\delta
(1)=e \ne 0$, we have $\dim(\ker (\delta ))\le \aut(\cB_{c-1})-1$.
On the other hand, using the hypothesis of induction to bound
$\aut(\cB_{c-1})$, we obtain

\begin{equation}\label{aut2}
\hom(\cB_{c-1},\cB_{c})=\hom(\cB_{c-1},\cO(a_{t+c-2}))+\dim(\im
(\eta))=\end{equation}
$$\hom(\cB_{c-1},\cO(a_{t+c-2}))+\dim(\ker(\delta))\le $$
$$\hom(\cB_{c-1},\cO(a_{t+c-2}))+\aut(\cB_{c-1})-1\le $$
$$\hom(\cB_{c-1},\cO(a_{t+c-2}))+1-1+\binom{h_0}{n+c}+ $$ $$
\sum_{i=1}^{c-4} \left (\sum _{r+s=i} \sum _{0\le i_1<...<i_r\le
t+i \atop 1\le j_1\le ... \le j_s\le t} (-1)^{i-r}
\binom{h_i+a_{i_1}+\cdots +a_{i_r}+b_{j_1}\cdots +b_{j_s} }{n+c}
\right )$$

\vskip 2mm \noindent where we set $h_i:=
2a_{t+1+i}+a_{t+2+i}+\cdots +a_{t+c-2}-\ell +n+c,$ for all
$i=1,\cdots , c-3$. Now, we will compute
$\hom(\cB_{c-1},\cO(a_{t+c-2}))$. To this end, we first observe
that $\Hom(\cB_{c-1},\cO)$ is the first Buchsbaum-Rim module
associated to

$$ \varphi_{c-1}^{*}: G_{c-1}^*= \bigoplus_{j=0}^{t+c-3}R(-a_j)
 \ra F ^{*}:=\bigoplus_{i=1}^tR(-b_i).$$
Therefore, we have the following free graded $R$-resolution

$$0 \ra \wedge^{t+c-2}G_{c-1}^* \otimes S_{c-3}(F)\otimes \wedge^t
F \ra \cdots \ra \wedge^{t+i+1} G_{c-1}^* \otimes S _{i}(F
)\otimes \wedge ^t F $$ $$ \ra \ldots \ra \wedge^{t+1}G_{c-1}^*
\otimes S_{0}(F)\otimes \wedge^t F  \ra \Hom(\cB_{c-1},\cO) \ra 0.
$$ Since, $$\wedge^t F =R(\sum_{i=1}^t b_i),$$ \vskip 2mm
$$S_m(F)=\bigoplus_{1\le j_1\le ...\le j_m\le
t}R(b_{j_1}+...+b_{j_m}), \mbox{ and}$$
\newpage $$ \wedge^r(G_{c-1}^*) = \bigoplus_{0\le i_1< ...< i_r\le
t+c-3}R(-a_{i_1}-...-a_{i_r}) \hspace{45mm}$$  $$ =
\bigoplus_{0\le i_1< ...< i_{t+c-2-r}\le
t+c-3}R(-\sum_{j=0}^{t+c-3}a_j+a_{i_1}+...+a_{i_{t+c-2-r}})  $$
$$\hspace{10mm} = \bigoplus_{0\le i_1< ...< i_{t+c-2-r}\le
t+c-3}R(-\sum_{j=0}^{t+c-2}a_j+a_{t+c+2}+a_{i_1}+...+a_{i_{t+c-2-r}})
$$ \noindent we have

$$ \wedge^{t+i+1}(G_{c-1}^*)\otimes S_i(F)\otimes \wedge ^tF  =
\hspace{120mm}$$ $$ \hspace{20mm} \bigoplus_{0\le i_1< ...<
i_{c-3-i}\le t+c-3 \atop 1\le j_1\le...\le j_i\le t
}R(-\ell+a_{t+c-2}+a_{i_1}+...+a_{i_{c-3-i}}+b_{j_1}+...+b_{j_i}).$$
 \vskip 4mm \noindent  So,  $ \dim _k ( \wedge^{t+i+1}(G_{c-1}^*)\otimes
S_i(F)\otimes \wedge ^tF) = $ $$ \hspace{30mm}\sum_{0\le i_1< ...<
i_{c-3-i}\le t+c-3 \atop 1\le j_1\le...\le j_i \le t}
\binom{-\ell+a_{t+c-2}+a_{i_1}+...+a_{i_{c-3-i}}+b_{j_1}+...+b_{j_i}+n+c}{n+c}
$$ \noindent and, we conclude that

\begin{equation} \label{comparar} \ hom(\cB_{c-1},\cO(a_{t+c-2}))= \hspace{50mm} \end{equation}
 $$\sum_{r+s=c-3} \left (\sum _{0\le
i_1< ...< i_{r}\le t+c-3 \atop 1\le j_1\le...\le j_s \le t}
(-1)^{c-3-r}\binom{h_{c-3}+a_{i_1}+...+a_{i_{r}}+b_{j_1}+...
+b_{j_s}}{n+c}\right )$$

\vskip 2mm \noindent being $h_{c-3}=2a_{t+c-2} -\ell +n+c$.
Putting altogether (\ref{aut1}), (\ref{aut2}) and
(\ref{comparar}), we obtain $$\aut(\cB_c)\le \binom{h_0}{n+c} +1+
$$ $$ \sum_{i=1}^{c-3}\left (\sum _{r+s=i} \sum _{0\le i_1< ...<
i_{r}\le t+i \atop 1\le j_1\le...\le j_s \le t } (-1)^{i-r}
\binom{h_i+a_{i_1}+\cdots +a_{i_r}+b_{j_1}\cdots +b_{j_s}
}{n+c}\right )$$

\vskip 2mm \noindent where we set $h_i:=
2a_{t+1+i}+a_{t+2+i}+\cdots +a_{t+c-2}-\ell +n+c$, for
$i=0,1,...,c-3$, which proves what we want.
\end{proof}

\begin{remark}\label{aut(B)=1}
Note that $\aut({\cB})=1$ provided $\ell
>2a_{t+c-2}+a_{t+c-3}+...+a_{t+1}+a_t$. (Indeed all binomials in
the expression in Proposition \ref{aut(B)} vanish.)
\end{remark}

 We are now ready to state the main result of this section.

\begin{theorem} \label{upperbound} Assume $a_0\le a_1\le ... \le a_{t+c-2}$,
$b_1\le ... \le b_t$ and $c\ge 2$. Set $\ell
:=\sum_{j=0}^{t+c-2}a_j-\sum_{i=1}^tb_i$ and $h_i:=
2a_{t+1+i}+a_{t+2+i}+\cdots +a_{t+c-2}-\ell +n+c$, for
$i=0,1,...,c-3$. It holds

\noindent (i) If $(c-1)a_{t+c-2}<\ell $, then
\[
 \dim W(\underline{b};\underline{a}) \le \sum_{i,j}
\binom{a_i-b_j+n+c}{n+c}- \sum _{i,j} \binom{a_i-a_j+n+c}{n+c}-
\]
\[
 \sum _{i,j} \binom{b_i-b_j+n+c}{n+c} + \sum _{j,i}
\binom{b_j-a_i+n+c}{n+c}+1.
\]

\noindent (ii) If $(c-1)a_{t+c-2}\ge \ell $, then
 \[
 \dim W(\underline{b};\underline{a}) \le \sum_{i,j}
\binom{a_i-b_j+n+c}{n+c}- \sum _{i,j} \binom{a_i-a_j+n+c}{n+c}-
\]
\[
 \sum _{i,j} \binom{b_i-b_j+n+c}{n+c} + \sum _{j,i}
\binom{b_j-a_i+n+c}{n+c}+\binom{h_0}{n+c} + 1+ \]
\[ \sum_{i=1}^{c-3}\left (\sum
_{r+s=i} \sum _{0\le i_1< ...< i_{r}\le t+i \atop 1\le
j_1\le...\le j_s \le t } (-1)^{i-r} \binom{h_i+a_{i_1}+\cdots
+a_{i_r}+b_{j_1}\cdots +b_{j_s} }{n+c}\right  ).
\]

\end{theorem}

\begin{proof}
 It follows from the inequality (\ref{boundw})
  and Proposition \ref{aut(B)}.
\end{proof}

\begin{remark}
{\rm Note that if $\ell
>2a_{t+c-2}+a_{t+c-3}+...+a_{t+1}+a_t$ then
 \[
 \dim W(\underline{b};\underline{a}) \le \sum_{i,j}
\binom{a_i-b_j+n+c}{n+c}- \sum _{i,j} \binom{a_i-a_j+n+c}{n+c}-
\]
\[
 \sum _{i,j} \binom{b_i-b_j+n+c}{n+c} + \sum _{j,i}
\binom{b_j-a_i+n+c}{n+c}+1.
\]
Indeed it follows form Theorem \ref{upperbound} and Remark
\ref{aut(B)=1}}.
\end{remark}

\begin{remark}
{\rm Given integers $a_0, a_1, \cdots , a_{t+c-2}$ and $b_1,
\cdots ,b_t$, we always have $\dim
W_s(\underline{b};\underline{a})=\dim
W(\underline{b};\underline{a})$. In fact, it
 is an easy consequence of Corollary \ref{WWs} and
the fact that an standard determinantal scheme is good
determinantal if it is generic a complete intersection and being a
generic complete intersection is an open condition.}
\end{remark}

\begin{example}(i) According to Ellingsrud Theorem (\cite{elli};
Theoreme 2), in codimension 2 case, the bound given in Theorem
\ref{upperbound} is sharp provided $n\ge 1$.

(ii) According to \cite{KMMNP}; Proposition 1.12, in codimension 3
case, the bound given in Theorem \ref{upperbound} is sharp,
provided $n\ge 1$ and $\depth _{I(Z_2)}D_2\ge 4$.

(iii) A {\em Rational normal Scroll} $X\subset \PP^N$ is a
non-degenerate variety of minimal degree (i.e.
$deg(X)=\codim(X)+1$) defined by the maximal minors of a $2\times
(c+1)$ matrix with linear entries ($c=\codim(X)$). As example of
rational normal scrolls we have the smooth, rational normal curves
of degree $d$ in $\PP^{d}$. It is well known that the family of
rational normal scrolls of degree $d$ and codimension $d-1$ in
$\PP^N$ is irreducible of dimension  $d(2n+2-d)-3$. So, again in
this case the bound given in Theorem \ref{upperbound} is sharp.

\end{example}

\vskip 2mm

We are led to pose the following questions.

\vskip 4mm

\begin{question} (i) When is the closure
of $W(\underline{b};\underline{a})$ an irreducible component of
$\Hi   ^{p(s)}(\proj{n+c})$?

(ii) Is $W(\underline{b};\underline{a})$ smooth or, at least,
generically smooth?

(iii) Under which extra assumptions are the bounds given in
Theorem \ref{upperbound} sharp?
\end{question}

We will address  questions (i) and (ii)  in next section 5 and
question (iii) in section 4.

\vskip 4mm Finally we want to point out that the inequality for
$\aut(\cB )$ in Proposition \ref{aut(B)} is indeed an equality.
One may see it by chasing the proof of Proposition \ref{aut(B)}
more carefully. We will, however, take the opportunity to compute
$\dim Aut(\cB )$ by a different method, leading to an apparently
new formula, and then prove that they coincide provided we assume
$a_0\le a_1 \le ... \le a_{t+c-2}$ and $b_1\le ... \le b_t$. This
new formula will be used in next section.

\begin{lemma} \label{222} With the previous notation, there is an exact sequence
$$0\rightarrow \Hom_R(B_c,F)\rightarrow \Hom_R(B_c,G_c)\rightarrow
\Hom_R(B_c,B_c)\rightarrow \Hom_R(M_c,M_c)\rightarrow 0.$$
\end{lemma}
\begin{proof} We apply $\Hom_R(B_c,.)$ to
$$0 \longrightarrow F\longrightarrow G_c\longrightarrow B_c
\longrightarrow 0 $$ and we get $$ 0\rightarrow
\Hom(B_c,F)\rightarrow \Hom_R(B_c,G_c)\rightarrow
\Hom_R(B_c,B_c)\rightarrow  $$ $$\rightarrow
\Ext^1_R(B_c,F)=M_c\otimes _R F \rightarrow
\Ext^1(B_c,G_c)=M_c\otimes _R G_c.$$ By
 Lemma \ref{key}, we have the exact sequence
$$0\rightarrow \Hom_R(M_c,M_c)\rightarrow M_c\otimes _R
F\rightarrow M_c\otimes _R G_c$$ and we are done.
\end{proof}

 Using the Buchsbaum-Rim
resolution $$0 \ra \wedge^{t+c-1}G_{c}^* \otimes S_{c-2}(F)\otimes
\wedge^t F \ra \cdots \ra \wedge^{t+i+1} G_{c}^* \otimes S _{i}(F
)\otimes \wedge ^t F $$ $$ \ra \ldots \ra \wedge^{t+1}G_{c}^*
\otimes S_{0}(F)\otimes \wedge^t F  \ra \Hom(B_{c},R) \ra 0, $$ we
immediately get the following Corollary:

\begin{corollary} Set $\tau _\nu := \hom_R(B_c,R)_{\nu }$.
Then,
$$\aut(B_c)=1+\sum_{j=0}^{t+c-2}\tau_{a_j}-\sum_{i=1}^t\tau_{b_i}.$$

\end{corollary}
\begin{proof} It follows from Lemmas \ref{key2} and  \ref{222} and the isomorphisms
$$\Hom_R(B_c,F)\cong \Hom(B_c,R)\otimes F\cong \oplus
_{i=1}^t\Hom(B_c,R(b_i)).$$
\end{proof}

\begin{proposition}
\label{2Aut(B)} Set $K_i:=$ $\hom(B_{i-1},R(a_{t+i-2}))_0$, for
$3\le i \le c$.  Suppose $b_1\le ... \le b_t$ and $a_0\le a_1 \le
... \le a_{t+c-2}$. Then we have $$\aut(\cB)=
1+K_3+K_4+\cdots+K_c,$$ and the inequality of $\dim Aut(\cB)$ in
Proposition \ref{aut(B)} is an equality.
\end{proposition}

\begin{proof}
Dualizing the exact sequence $0\rightarrow R(a_{t+c-2})\rightarrow
B_c\rightarrow B_{c-1}\rightarrow 0$, we get $$ 0\rightarrow
\Hom(B_{c-1},R)\rightarrow \Hom_R(B_c,R)\rightarrow
R(-a_{t+c-2})\rightarrow  M_{c-1} \rightarrow M_c\rightarrow 0$$
which together with (\ref{Mi}) gives us the exact sequence
\begin{equation}\label{00} 0\rightarrow \Hom(B_{c-1},R)\rightarrow
\Hom_R(B_c,R)\rightarrow I_{D_{c-1}}(-a_{t+c-2})\rightarrow
0.\end{equation}

Look at the commutative diagram
\[
\begin{array}{cccccc}
0 & & 0 \\\downarrow  & & \downarrow  &     \\
 \Hom(B_{c-1},F)_0 &  \longrightarrow &   \Hom(B_{c-1},G_c)_0
\\ \downarrow  & &
\downarrow  &     \\ \Hom(B_{c},F)_0 &  \longrightarrow &
\Hom(B_{c},G_c)_0
\\ \downarrow  & & \downarrow   \\
0=\Hom(R(a_{t+c-2}),F)_0 &  \longrightarrow &
\Hom(R(a_{t+c-2}),G_c)_0.

\end{array}
\]
By (\ref{00}), $\Hom(B_{c},G_c)_0\rightarrow
\Hom(R(a_{t+c-2}),G_c)_0$ is zero because its image is
$(I_{D_{c-1}}\otimes G_c(-a_{t+c-2}))_0$. Hence, we get
$$\hom(B_c,G_c)_0-\hom(B_c,F)_0=$$
$$\hom(B_{c-1},G_{c-1})_0+\hom(B_{c-1},R(a_{t+c-2}))_0-\hom(B_{c-1},F)_0.$$
Since, we have $$\aut(\cB_c)=1+\hom(B_c,G_c)_0-\hom(B_c,F)_0$$ by
Lemma \ref{222} and we may suppose
$$\aut(\cB_{c-1})=1+\hom(B_{c-1},G_{c-1})_0-\hom(B_{c-1},F)_0$$ we
have proved
\begin{equation}\label{newnew}
\aut(\cB_c)=K_c+\aut(\cB_{c-1}).\end{equation}
 Now, we conclude by
induction taking into account that $\aut(\cB_2)=$ $ \hom
(I_{D_2},I_{D_2})_0=1$. Moreover combining (\ref{newnew}) and the
definition of $K_c$ with  (\ref{comparar}) we see that the
expression we have got for $\aut (\cB_c)$ coincides with the
corresponding binomials in the expression of $\aut (\cB _c)$ in
Proposition \ref{aut(B)}, and it follows that the inequality must
be an equality.
\end{proof}

So, we can rewrite Theorem \ref{upperbound} and we have

\begin{proposition}\label{bound2} With the above notation
\[
 \dim W(\underline{b};\underline{a}) \le \sum_{i,j}
\binom{a_i-b_j+n+c}{n+c}- \sum _{i,j} \binom{a_i-a_j+n+c}{n+c}-
\]
\[
 \sum _{i,j} \binom{b_i-b_j+n+c}{n+c} + \sum _{j,i}
\binom{b_j-a_i+n+c}{n+c}+1+K_3+\cdots +K_c.
\]
\end{proposition}

\begin{proof} It follows from the inequality (\ref{boundw})
  and Proposition \ref{2Aut(B)}.
\end{proof}

\begin{remark} {\rm One may show that the right hand side of the
inequality for $\dim W(\underline{b};\underline{a})$ in
Proposition \ref{bound2} is equal to dim $ \Ext^1_R(B_c,B_c)_0$.
This indicates an interesting connection to the deformations of
the $R$-module $B_c$.}
\end{remark}


\section{The dimension of the determinantal locus}

The purpose of this section is to analyze when the bound given in
Theorem \ref{upperbound} is sharp. We will see that under mild
conditions the upper bound of dim$W(\underline{b};\underline{a})$
given in the preceding section is indeed the dimension of the
determinantal locus $W(\underline{b};\underline{a})$ provided the
codimension $c$ is small. Indeed, if $2\le c \le 3$ and $n\ge 1$,
this is known (\cite{KMMNP}, \cite{elli}) while for $4\le c \le 5$
is a consequence of the main theorem of this section. If $c\ge 6$
we also get the expected dimension formula for
$W(\underline{b};\underline{a})$ under more restrictive
assumptions. As in the preceding section the proofs use induction
on $c$ by successively deleting columns of the largest possible
degree.

\vskip 2mm We keep the notation introduced in \S 2 and \S 3: see
in particular (\ref{gradedmorfismo})-(\ref{DiMi}). If we denote by
$W(F,G):=W(\underline{b};\underline{a})$ and by $\VV
(F,G_i):=\Hom_{{\cO}_{\PP^{n+c}}}(\tilde{F},\tilde{G_i})$ the
affine scheme whose rational points are the morphisms  from
$\tilde{F}$ to $\tilde{G_i}$, we have by the definition of $W
(F,G_c)$ and $W (F,G_{c-1})$ a diagram of rational maps

 $$
\xymatrix{
 \VV(F,G_c)\ar[r] \ar@{->>}[d] & \VV
(F,G_{c-1})\ar@{->>}[d]
\\ W(F,G_c)& W(F,G_{c-1})
}
 $$

\noindent where the down arrows are dominating and rational and
$\VV(F,G_c) \longrightarrow  \VV(F,G_{c-1})$ is defined by
deleting the last column.

To prove that the upper bound of dim$W(F,G_c)$ of Proposition
\ref{bound2} is also a lower bound, we need a
deformation-theoretic technical result which computes the
dimension of $W(F,G_c)$ in terms of the dimension of
$W(F,G_{c-1})$. To do so, we consider the Hilbert flag scheme
$D(p,q)$ parameterizing "pairs" $X\subset Y$ of closed subschemes
of $\PP^{n+c}$ with Hilbert polynomial $p$ and $q$ respectively
and the subset $D(F,G_i,G_{i-1})$ of "pairs" $X\subset Y$  where
$X\in W(F,G_i)$ is a good determinantal scheme defined by a matrix
$\cA _i\in \VV(F,G_i)$ and $Y$ is a good determinantal scheme
defined by the matrix $\cA _{i-1}\in \VV(F,G_{i-1})$ obtained by
deleting the last column of $\cA _i$. Then the diagram above fits
into
 $$
\xymatrix{ \VV(F,G_c)\ar[rrr] \ar@{->>}[rd]\ar@{->>}[dd] &&& \VV
(F,G_{c-1}) \ar@{->>}[d]  \\ &D(F,G_c,G_{c-1})\ar@{->}[ld]^{p_1}
\ar@{->}[rr]_{p_2}&&W(F,G_{c-1})\\
W(\underline{b};\underline{a})=W(F,G_c)&&&
\\ }
 $$

\noindent where $p_1$ and $p_2$ are the restriction of the natural
projections $pr_1:D(p,q)\longrightarrow \Hi ^p(\PP^{n+c})$ and
$pr_2:D(p,q)\longrightarrow \Hi ^q( \PP^{n+c})$ respectively, and
where $\VV(F,G_c)\twoheadrightarrow D(F,G_c,G_{c-1})$ is
dominating and rational by definition. Denoting $$m_i(\nu )=\dim
_kM_i(a_{t+i-2})_{\nu}$$ we have

\begin{proposition}\label{main0} Let $c\ge 3$.
Suppose that $W(\underline{b};\underline{a})\ne \emptyset$ and
that $\depth_{I(Z_{c-1})}D_{c-1}\ge 2$ for a general enough
$D_{c-1}\in W(F;G_{c-1})$. Then
\begin{itemize}
\item[(1)] $p_2$ is dominating and
$$\dim D(F,G_c,G_{c-1})\ge \dim W(F,G_{c-1})+m_c(0);$$
\item[(2)] $\dim W(\underline{b};\underline{a})\ge
\dim
W(F,G_{c-1})+m_c(0)-\hom_{\cO_{\PP^{n+c}}}(\cI_{X_{c-1}},\cI_{c-1})$.

\end{itemize}
\end{proposition}
\begin{proof}

Due to \cite{KMNP}, Proposition 3.2, we see that for any $Y\in
W(F,G_{c-1})$ there exists a regular section
$R/I(Y)\hookrightarrow M_{c-1}(a_{t+c-2})$ whose cokernel is
supported at some $X$ with $\dim X < \dim Y$, and such that
$M_{c-1}$ is the cokernel of the morphism $\varphi
^*_{c-1}:G^*_{c-1}\rightarrow F^*$ as in \S 3. Moreover, for a
given $Y$, the mapping cone construction shows that for any
regular section $R/I(Y)\hookrightarrow M_{c-1}(a_{t+c-2})$ there
is a morphism $\varphi ^*_{c}:G^*_{c}\rightarrow F^*$ which
reduces to the given $\varphi ^*_{c-1}$ by deleting the extra (say
the last) column of the corresponding matrix. This shows that
$p_2$ is dominating and that the fibers $p_2^{-1}(Y)$ "contains"
the space of regular sections of $M_{c-1}(a_{t+c-2})$ in a natural
way.

\vskip 2mm More precisely note that any $Y\in W(F,G_{c-1})$ has
the same Betti numbers of a Buchsbaum-Rim resolution (cf.
Proposition 2.2) and so the same dimension of
$M_{c}(a_{t+c-2})_0=M_{c-1}(a_{t+c-2})_0/k$. If $Y$ is general
enough, we have $$ \dim D(F,G_c,G_ {c-1})=\dim W(F;G_{c-1})+\dim
p_2^{-1}(Y)$$ by generic flatness. Hence
 it suffices to see that $\dim p_2^{-1}(Y)\ge m_c(0)$. Pick
 $(X\subset Y)\in p_2^{-1}(Y)$, look at (\ref{DiMi}) and consider
 the injection $M_c(a_{t+c-2})_0\hookrightarrow (N_{X/Y})_0$. In
 the tangent space $(N_{X/Y})_0$ of $pr_2^{-1}(Y)\supseteq
 p_2^{-1}(Y)$ at $(X\subset Y)$ we therefore have a
 $m_c(0)$-dimensional family arising from deforming the matrix
 $\cA =[\cA _{c-1},L]$ of $\varphi _c^*$ leaving $\varphi
 _{c-1}^*$ (i.e. $\cA _{c-1}$) fixed ($L$ is the last column of
 $\cA $). We may think of the last column of such a deformation of
 $\varphi_c^*$ as $L+\sum_{i=1}^{m_c(0)}t_iL^{(i)}$ mod. $(t_1,t_2,...,t_{m_c(0)})^2$
 where the $t_i's$ are indeterminates and where the degree
 matrix of the columns $L^{(i)}$ are exactly the same as that of
 $L$. Since the degeneracy locus of the $t\times (t+c-1)$ matrix
 $[\cA _{c-1},L+\sum_{i=1}^{m_c(0)}t_iL^{(i)}]$ defines a flat
 family over some open subset $T$ of $Spec(k[t_1,...,t_{m_c(0)}])$
 containing the origin (because the Eagon-Northcott complex over
 $Spec(k[\underline{t}])$ must be acyclic over some $T$ provided
 the pullback to $(0)\in Spec(k[\underline{t}])$ is acyclic), we
 see that the  fiber $p_2^{-1}(Y)$ contains a $m_c(0)$-dimensional
 (linear) family, as required. This proves (1).

\vskip 4mm
 (2) It is straightforward to get (2) from (1). Indeed,
 $$\dim D(F,G_c,G_{c-1})-\dim W(\underline{b};\underline{a})\le
 \dim p_1^{-1}(D_c)$$
 and since $p_1^{-1}(D_c)$ is contained in the full fiber of the
 first projection $pr_1:D(p,q)\rightarrow\Hi ^p(\PP^{n+c})$ whose
 fiber dimension is known to have $\hom(\cI_{X_{c-1}},\cI_{c-1})$
  as an upper bound (e.g. \cite{KMMNP}, Chapter
 9), we easily conclude.
\end{proof}

Proposition \ref{main0} allows us, under some assumptions, to find
a lower bound for $\dim W(\underline{b};\underline{a})$ provided
we have a lower bound of $\dim W(F,G_{c-1})$. Indeed, since it is
easy to find $m_{i}(0)$ using the Buchsbaum-Rim resolution of
$M_i$ or by using (\ref{Mi}) recursively, it remains to find $\hom
(\cI_{X_{i}},\cI_{i})$ in terms of $\hom (\cI_{X_{i-1}},\cI
 _{i-1})$.

\begin{lemma} \label{ext}
Set $a=a_{t+i-2}-a_{t+i-1}$. \begin{itemize} \item[(a)]
If
$\Ext^1_{D_{i-1}}(I_{D_{i-1}}\otimes I_{i-1}^*,I_{i-1})_{\nu
+a}=0$ and $\depth_{I(Z_{i-1})}D_{i-1}\ge 3$, then
$$\hom(I_{D_{i}}, I_{i})_\nu \le \dim
(D_i)_{\nu+a}+\hom(I_{D_{i-1}}, I_{i-1 })_{\nu +a}.$$
 \item[(b)]
If $\Ext^2_{D_{i-1}}(I_{D_{i-1}}\otimes I_{i-1}^*,I_{i-1})_{\nu
+a}=0$ and $\depth_{I(Z_{i-1})}D_{i-1}\ge 4$, then $$\Ext^1
_{D_{i-1}}(I_{D_{i-1}}/I^2_{D_{i-1}},I_{i-1})_{\nu
+a}=0\Rightarrow  \Ext^1
_{D_{i}}(I_{D_{i}}/I^2_{D_{i}},I_{i})_{\nu }=0.$$
\end{itemize}

\end{lemma}
\begin{remark} \label{rem} {\rm Since $I_{i-1}=M_{i-1}(a_{t+i-2})^*$,  we
have $$\Ext^1_{D_{i-1}}(I_{D_{i-1}}\otimes
I_{i-1}^*,I_{i-1}(a))\cong \Ext^1_{D_{i-1}}(I_{D_{i-1}}\otimes
M_{i-1},M_{i-1}^*(-a_{t+i-2}-a_{t+i-1})).$$}
\end{remark}

\begin{proof} (a) We consider the two exact sequences
\begin{equation}\label{1}
0\rightarrow \Hom _R(I_{i-1}, I_{i}) \rightarrow \Hom
_R(I_{D_{i}}, I_{i}) \rightarrow \Hom _R(I_{D_{i-1}},
I_{i})\end{equation}
\begin{equation}\label{2}
0\rightarrow \Hom _{D_{i-1}}(I_{D_{i-1}}\otimes I_{i-1}^*,I_{i-1})
\rightarrow \Hom _{D_{i-1}}(I_{D_{i-1}}\otimes
I_{i-1}^*,D_{i-1})\rightarrow \Hom _{D_{i-1}}(I_{D_{i-1}}\otimes
I_{i-1}^*,D_i).
\end{equation}
We have
$\depth_{I(Z_{i-1})}D_{i-1}\ge 3$ and hence
$\depth_{I(Z_{i-1})}D_{i}\ge 2$ and $\depth_{I(Z_{i-1})}I_{i}\ge
2$ and we get by (\ref{NM})

\begin{equation}
\label{new}\Hom(I_{i-1}, I_{i})\cong H^0_{*}(U_{i-1},\cH
om(\cI_{i-1},\cI_i))\cong \end{equation}
 $$H^0_{*}(U_{i-1},\cH
om_{\cO _{X_{i}}}(\cI _{i-1}\otimes_{\cO _{X_{i-1}}}\cO_
{X_i}\otimes \cI ^* _i, \cO _{X_{i}}))\cong D_i(a)$$ because, by
(\ref{Mi}), $\tilde{M}_{i-1}(a_{t+i-2})\otimes_{ \cO_{X_{i-1}}}
\cO_{X_{i}} |_{U_{i-1}} \cong
\tilde{M}_i(a_{t+i-1})(a)|_{U_{i-1}}$ and hence $$\cI
_{i-1}\otimes_{\cO _{X_{i-1}}}\cO_{X_{i}} |_{U_{i-1}}\cong
\cI_{i}(-a)|_{U_{i-1}}.$$

For similar reasons; $$\Hom_{D_{i-1}}(I_{D_{i-1}}\otimes
I_{i-1}^*,D_{i-1})\cong H^0_{*}(U_{i-1},\cH om(\cI _{X_{i-1}}/\cI
^2_{X_{i-1}}\otimes \cI^*_{i-1}, \cO_{X_{i-1}}))\cong $$ $$
H^0_{*}(U_{i-1},\cH om(\cI _{X_{i-1}}/\cI ^2_{X_{i-1}}, \cI_{i-1})
)\cong \Hom _R(I_{D_{i-1}},I_{i-1})$$
 and
$$\Hom_R(I_{D_{i-1}},I_{i})\cong H^0_{*}(U_{i-1},\cH
om_{\cO_{X_{i}}}(\cI _{X_{i-1}}/\cI ^2_{X_{i-1}}\otimes \cI^*_{i},
\cO_{X_{i}}))$$ is further isomorphic to
$$\Hom_{D_{i-1}}(I_{D_{i-1}}\otimes I_{i-1}^*,D_i(a))\cong
H^0_{*}(U_{i-1}, \cH om_{\cO_{X_{i}}}(\cI _{X_{i-1}}/\cI
^2_{X_{i-1}}\otimes \cI^*_{i-1}\otimes _{\cO_{X_{i-1}}}
\cO_{X_{i}}, \cO_{X_{i}}(a))).$$

Putting all this together, we get that the  exact sequences
(\ref{1}) and (\ref{2}) reduce, in degree $\nu$ and $\nu +a$
resp., to

\begin{equation} \label{47} 0\rightarrow (D_i)_{\nu +a}\rightarrow
\Hom_R(I_{D_i},I_i)_\nu \rightarrow \Hom _R(I_{D_{i-1}},I_i)_{\nu
} \cong \Hom_{D_{i-1}}(I_{D_{i-1}}\otimes
I_{i-1}^*,D_i)_{\nu+a}\rightarrow
\end{equation}
\begin{equation}\label{48}0 \rightarrow \Hom(I_{D_{i-1}}
\otimes I_{i-1}^*,I_{i-1} )_{\nu+a}\rightarrow
\Hom_R(I_{D_{i-1}},I_{i-1})_{\nu+a}\rightarrow
\Hom_{D_{i-1}}(I_{D_{i-1}}\otimes
I_{i-1}^*,D_i)_{\nu+a}\rightarrow 0
\end{equation}
where (\ref{48}) is short-exact by assumption.  Taking dimensions,
we immediately get (a).

\vskip 2mm (b) As in (\ref{new}) we see that
\begin{equation}\label{46new}
\Ext ^1_{D_{i-1}}(I_{i-1},I_i)\cong H^1_*(U_{i-1},\cO
_{X_i}(a))=0.
\end{equation}
Sheafifying (\ref{47}) and (\ref{48}) and taking global sections,
we get
{\scriptsize
\begin{equation}\label{47new}
0  \rightarrow  H^1_*(U_{i-1},\cH om(\cI_{X_{i}}/\cI
^2_{X_{i}},\cI_{i}(-a)))  \rightarrow  H^1_*(U_{i-1},\cH
om(\cI_{X_{i-1}}/\cI ^2_{X_{i-1}},\cI_{i}(-a)))   \rightarrow
\end{equation} $$ \hspace{40mm} \| $$ $$  \rightarrow  H^1_*(U_{i-1},\cH
om(\cI_{X_{i-1}}\otimes \cI^*_{i-1},\cO_{X_{i-1}})) \rightarrow
H^1_*(U_{i-1},\cH om(\cI_{X_{i-1}}\otimes \cI
^*_{i-1},\cO_{X_{i}}))  \rightarrow H^2_*(U_{i-1},\cH
om(\cI_{X_{i-1}}\otimes \cI ^*_{i-1},\cI_{i-1})). $$ }

\noindent Since $\cH om(\cI_{X_{i-1}}\otimes  \cI^*_{i-1},\cO
_{X_{i-1}})\cong \cH om(\cI_{X_{i-1}}/\cI
^2_{X_{i-1}},\cI_{i-1})$, then $\depth _{I(Z_{i-1})}D_{i-1} \ge 4$
and (\ref{NM}) show that the $H_*^{i}$-groups of (\ref{47new}) are
isomorphic to the $\Ext ^{i}$-groups in the following diagram
\begin{equation}\label{48new}\hspace{-20mm}
0\rightarrow \Ext^1_{D_i}(I_{D_i}/I^2_{D_i},I_{i})\rightarrow
\Ext^1_{D_{i-1}}(I_{D_{i-1}}/I^2_{D_{i-1}},I_{i})\rightarrow
\end{equation}
$$ \hspace{25mm} \| $$ $$ \rightarrow
\Ext^1_{D_{i-1}}(I_{D_{i-1}}/I^2_{D_{i-1}},I_{i-1}(a))\rightarrow
\Ext^1_{D_{i-1}}(I_{D_{i-1}}\otimes I^*_{i-1},D_{i}(a))\rightarrow
\Ext^2_{D_{i-1}}(I_{D_{i-1}}\otimes I^*_{i-1},I_{i-1}(a))$$ of
exact horizontal sequences. Using (\ref{48new}) we easily get (b).
\end{proof}

\begin{remark} \label{rem44}
{\rm By (\ref{47}) the conclusion of Lemma 4.2 obviously holds provided
 we have $\Hom_R(I_{D_{i-1}},I_{i})_{\nu}=0$. Using the Eagon-Northcott
resolution of $I_{D_{i-1}}$ (i.e. of $D_{i-1}$), one may see that
this $\Hom_\nu$-group vanishes if $a_{t+i-2}$ is large enough.}
\end{remark}

\vskip 4mm Put $$\lambda_c:= \sum_{i,j}  \binom{a_i-b_j+n+c}{n+c}
+ \sum_{i,j} \binom{b_j-a_i+n+c}{n+c} - $$ $$ \sum _{i,j}
\binom{a_i-a_j+n+c}{n+c}-
 \sum _{i,j} \binom{b_i-b_j+n+c}{n+c}  + 1$$ where the indices  belonging to
$a_j$ (resp. $b_i$) ranges over $0\le j \le t+c-2$ (resp. $1\le i
\le t$). We define $\lambda _{c-1}$ by the analogous expression
where now $a_j$ (resp $b_i$) ranges over $0\le j \le t+c-3$ (resp.
$1\le i \le t$). It follows after a straightforward computation
that
\begin{equation}\label{lambdas}
\lambda _c= \end{equation} $$\lambda_{c-1}+\sum_{i=1}^t
\binom{a_{t+c-2}-b_i+n+c}{n+c}- \sum _{j=0}^{t+c-3}
\binom{a_{t+c-2}-a_j+n+c}{n+c}-
 \sum _{j=0}^{t+c-2} \binom{a_j-a_{t+c-2}+n+c}{n+c}.
$$

We now come to the main theorem of this section which shows that
the inequalities in Theorem \ref{upperbound} are equalities under
certain assumptions. Recalling the equivalent expression of the
upper bound of $\dim W(\underline{b};\underline{a})$ given in
Proposition \ref{bound2}, we have

\begin{theorem}\label{MAIN}
Let $a_0\le a_1\le ... \le a_{t+c-2}$ and $b_1\le...\le b_t$ and
assume $a_{i-min(2,t)}\ge b_{i}$ for $min(2,t)\le i \le t$. Let
$c\ge 3$ and let $W(\underline{b};\underline{a})$ be the locus of
good determinantal schemes in $\PP^{n+c}$ where $n\ge 0$ if $c\ge
4$ and $n\ge 1$ if $c=3$. For a general $Proj(A)\in
W(\underline{b};\underline{a})$, let $R\twoheadrightarrow
D_2\twoheadrightarrow D_3\twoheadrightarrow ... \twoheadrightarrow
D_c=A$ be the flag obtained by successively deleting columns from
the right hand side. If $$\Ext ^1_{D_{i-1}}(I_{D_{i-1}}\otimes
I_{i-1}^*,I_{i-1})_{\nu }=0 \mbox{ for } \nu \le 0 \mbox{ and }
3\le i\le c-1$$ then $$ \dim W(\underline{b};\underline{a})=
\lambda _c+K_3+K_4+...+K_c$$ where $K_i=\hom
(B_{i-1},R(a_{t+i-2}))_0$ for $3\le i\le c$.
\end{theorem}

\begin{remark} If $c=2$ and $n\ge 1$ one knows by \cite{elli} that $$ \dim
W(\underline{b};\underline{a})= \lambda _2.$$ The same formula
holds if $c=2$ and $n=0$ as well. In this case one may get the
formula by taking a general $Proj(A)\in
W(\underline{b};\underline{a})$ and show that $$\hom _R(I_A,A)_0=
\ext^1_R(I_A,I_A)_0=\lambda _2$$ by e.g. using \cite{KP}, (26). We
leave the details to the reader.
\end{remark}
\begin{proof} Due to Remark \ref{dep} and the assumption
$a_{i-min(2,t)}\ge b_{i}$ for $min(2,t)\le i \le t$, the set
$Z_i=Sing(X_i)$ satisfies $\depth_{I(Z_i)}D_i\ge 3$ for $2\le i
\le c-2$, $\depth_{I(Z_{c-1})}D_{c-1}\ge 2$ (and also
$\depth_{I(Z_2)}D_2\ge 3$ in case $c=3$ since $n\ge 1$) by
choosing $C=Proj(A)$ general enough in $
W(\underline{b};\underline{a}).$

To use Proposition \ref{main0}, we only need to compute $m_c(0)$
and $\hom(I_{D_{c-1}},I_{c-1})_0$ because we may by induction
suppose that $\dim W(F,G_{c-1})=\lambda _{c-1}+K_3+...+K_{c-1}$
for $c\ge 3$ (interpreting the case $c-1=2$ as $\lambda _2$).  By
(\ref{Mi}) and (\ref{defMi}) we get
\begin{equation}\label{Mc0}
m_0(c)=\dim M_{c-1}(a_{t+c-2})_0-1=
\end{equation}
$$ \dim F^*(a_{t+c-2})_0-\dim
G^*_{c-1}(a_{t+c-2})_0+\hom(B_{c-1},R(a_{t+c-2}))_0-1=$$
$$=\sum_{i=1}^t \binom{a_{t+c-2}-b_i+n+c}{n+c}- \sum
_{j=0}^{t+c-3} \binom{a_{t+c-2}-a_j+n+c}{n+c}+K_c-1.$$

Thanks to Lemma \ref{ext}, we can find an upper bound of
$\hom(I_{D_{c-1}},I_{c-1})_0$. We have
$$\hom(I_{D_{c-1}},I_{c-1})_0\le
\binom{a+n+c}{n+c}+\hom(I_{D_{c-2}},I_{c-2})_a$$ because
$a=a_{t+c-3}-a_{t+c-2}\le 0$ and $\dim(D_i)_a$, which is either 0
or 1, must be equal to the binomial coefficient above. Repeatedly
use of Lemma \ref{ext} implies
\begin{equation}\label{4111}\hom(I_{D_{c-1}},I_{c-1})_0\le
 \sum _{j=t+1}^{t+c-3}
\binom{a_j-a_{t+c-2}+n+c}{n+c}+
\hom(I_{D_{2}},I_{2})_{a_{t+1}-a_{t+c-2}}.\end{equation}

It remains to compute $\hom(I_{D_{2}},I_{2})_{\alpha }$ with
$\alpha =a_{t+1}-a_{t+c-2}$. Using (\ref{NM})  (cf. the proof of
Lemma 4.2), we get $$\Hom(I_{D_{2}},I_{2})\cong
\Hom_{D_2}(I_{D_{2}}\otimes I_{2}^*,D_2)\cong
\Hom_{D_2}(I_{D_{2}}\otimes M_2(a_{t+1}),D_{2}).$$

\noindent Moreover, if $\ell _2=\sum _{j=0}^t a_j-\sum _{i=1}^t
b_i$, then $M_2\cong K_{D_2}(-\ell _2+n+c+1)$ by Proposition 2.2.
In codimension $c=2$, one knows $$(I_{D_2}/I_{D_2}^2)^*\cong
\Ext^1_R(I_{D_2},I_{D_2})\cong\Ext^1_R(I_{D_2},D_2)\otimes
I_{D_2}\cong K_{D_2}(n+c+1)\otimes I_{D_2}$$ and since $\depth
_{I(Z_2)}D_2 \ge 3$ and hence  $\depth
_{I(Z_2)}I_{D_2}/I_{D_2}^2\ge 2$ (because the codepth of $
I_{D_2}/I_{D_2}^2$ is $\le 1$ by \cite{AH}), we get

\begin{equation} \label{412} \Hom(I_{D_2},I_2)_{\alpha} \cong \Hom
(I_{D_2}\otimes K_{D_2}(n+c+1),D_2)(\ell _2-a_{t+1})_{\alpha }
\cong
\end{equation}
$$\cong (I_{D_2}/I_{D_2}^2)^{**}(\ell _2-a_{t+1})_{\alpha } \cong
(I_{D_2}/I_{D_2}^2)_{\ell _2-a_{t+c-2} }.$$
 Thus the inequality $a_j \le a_{t+c-2}$ and the exact
sequences

$$0\rightarrow F\rightarrow G_2=\oplus _{j=0}^t R(a_j)\rightarrow
I_{D_2}(\ell _2)\rightarrow 0$$

\begin{equation} \label{2power}
0 \rightarrow\wedge ^2 F\rightarrow F\otimes G_2\rightarrow
S_2G_2\rightarrow I^2_{D_2}(2\ell _2)\rightarrow 0
\end{equation}
show $$\hom(I_{D_2},I_2)_{\alpha} =\dim (G_2)_{-a_{t+c-2}}=\sum
_{j=0}^t\binom{a_j-a_{t+c-2}+n+c}{n+c}.$$

Using this last inequality together with (\ref{Mc0}), (\ref{4111})
and Proposition 4.1, we get by induction $$ \dim
W(\underline{b};\underline{a})\ge \lambda
_{c-1}+K_3+...+K_{c-1}+\sum _{i=1}^t\binom{a_{t+c-2}-b_i+n+c}{n+c}
$$ $$ -\sum
_{j=0}^{t+c-3}\binom{a_{t+c-2}-a_j+n+c}{n+c}+K_c-1-\sum
_{j=0}^{t+c-3}\binom{a_j-a_{t+c-2}+n+c}{n+c}$$ $$=\lambda
_c+K_3+K_4+...+K_c$$ where the last equality is due to
(\ref{lambdas}). Combining with Proposition \ref{bound2}, we get
the Theorem.
\end{proof}
Note that the vanishing assumption of Theorem 4.5 is empty if
$c=3$. Hence, we have

\begin{corollary}\label{cod3}
 Let
$W(\underline{b};\underline{a})$ be the locus of good
determinantal schemes in $\PP^{n+c}$ where $n\ge 1$ and $c= 3$. If
$a_{i-min(2,t)}\ge b_{i}$ for $min(2,t)\le i \le t$, then $$ \dim
W(\underline{b};\underline{a})= \lambda _3+K_3.$$
\end{corollary}
 \qed

\begin{remark} {\rm The above Corollary essentially generalizes \cite{KMMNP},
Corollary 10.15(i) where the depth condition is slightly stronger
than the one we use in the proof of Theorem 4.5. The only missing
part is that the assumption $n\ge 1$ excludes the interesting case
of 0-dimensional good determinantal schemes. See Corollary
\ref{cod3dim0} for the 0-dimensional case.}
\end{remark}

\vskip 2mm To apply Theorem \ref{MAIN} in the codimension $c=4$
case, it suffices to prove that $$\Ext^1_{D_2}(I_{D_2}\otimes
I_2^*,I_2)=0.$$

\noindent Due to Remark \ref{rem} and Proposition 2.2, the
$\Ext^1$-group above is isomorphic to a twist of

\begin{equation}\label{new?}\Ext^1_{D_2}(I_{D_2}\otimes M_2,M_2^*)\cong
\Ext^1_{D_2}(I_{D_2}\otimes K_{D_2},K_{D_2}^*)(2\ell
_2-2n-2c-2).\end{equation}

\noindent Hence,  all we need follows from

\begin{lemma} \label{key4} Let $R\twoheadrightarrow D=R/I_D$ be a
Cohen-Macaulay codimension 2 quotient and suppose
$Proj(D)\hookrightarrow \PP^{n+c}$ is a local complete
intersection outside a closed subset $Z\subset Proj(D)$ which
satisfies $\depth_{I(Z)}D\ge 4$. Then $$\depth _{\goth m} \Hom
_D(I_D\otimes K_D,K_D^*)\ge \depth _{\goth m}D-1.$$

In particular,  $\depth _{I(Z)}\Hom _D(I_D\otimes K_D,K_D^*)\ge 3$
and hence $$\Ext^1_D(I_D\otimes K_D,K_D^*)=0.$$

\end{lemma}
\begin{proof} $D$ is determinantal, say $D=D_2$ and we have a
minimal free $R$-resolution
\begin{equation}\label{resol} 0\rightarrow F \rightarrow
G_2\rightarrow I_D(\ell _2)\rightarrow 0
\end{equation}
as previously. If $H_i$ is the $i$-th Koszul homology built on
some set of minimal generators of $I_D$, it suffices to show that
there are two exact sequences
\begin{equation}\label{415} 0\rightarrow
\Hom_D(K_D(n+c+1),H_1)\rightarrow\wedge^2 (F(-\ell _2))\otimes
D\rightarrow H_2\rightarrow 0
\end{equation}
\begin{equation} \label{416}
0\rightarrow \Hom_D(K_D,H_1)\rightarrow K_D^*\otimes G_2(-\ell
_2)\rightarrow \Hom(I_D\otimes K_D(n+c+1),K_D^*)\rightarrow 0.
\end{equation}

\noindent Indeed, $H_i$ are maximal Cohen-Macaulay modules by
\cite{hun}. Hence, the first sequence shows that
$\Hom_D(K_D(n+c+1),H_1)$ is maximal Cohen-Macaulay while the
second shows that the codepth of $\Hom(I_D\otimes
K_D(n+c+1),K_D^*)$ is at most 1 and all conclusions of  the lemma
follow easily (cf. (\ref{NM}) for the last conclusion).

To see that (\ref{415}) is exact we deduce, from (\ref{resol}),
the exact sequence $$ 0\rightarrow K_D(n+c+1)^*\rightarrow F(-\ell
_2)\otimes _R D\rightarrow G_2(-\ell _2)\otimes _R D\rightarrow
I_D/I_D^2\rightarrow 0.$$ Indeed, we only need to prove that $
K_D(n+c+1)^*=\ker [F(-\ell _2)\otimes _R D\rightarrow G_2(-\ell
_2)\otimes _R D]$ which follows by applying $\Hom_R(.,D)$ to $$
..\rightarrow G_2(-\ell _2)^*\rightarrow F(-\ell _2)^*\rightarrow
\Ext^1_R(I_{D},R)\cong K_D(n+c+1)\rightarrow 0.$$

\noindent Since one moreover knows
\begin{equation}\label{417}
0\rightarrow H_1\rightarrow G_2(-\ell _2)\otimes _R D\rightarrow
I_D/I^2_D\rightarrow 0
\end{equation}
we get the exact sequence
\begin{equation}\label{418}
0\rightarrow K_D(n+c+1)^*\rightarrow F(-\ell _2)\otimes _R
D\rightarrow H_1\rightarrow 0,
\end{equation}
from which  we see that the Cohen-Macaulayness of $K_D(n+c+1)^*$
follows from that of $H_1.$ Sheafifying (\ref{418}) and using
\cite{har}, Ch II, exec. 5.16, we get an exact sequence

$$0\rightarrow {\tilde K_D(n+c+1)}^*\otimes {\tilde
H_1}|_{U}\rightarrow \wedge ^2({\tilde  F(-\ell _2))}\otimes
\tilde{ D}|_{U}\rightarrow \wedge ^2{\tilde H_1}|_{U} $$ where
$U=Proj(D)-Z$. Applying $H_{*}^0(U,.)$ and recalling that
$H^0_{*}(U,\wedge ^2{\tilde H_1})\cong H_2$
 \cite{KP}, Proposition 18, we get the exact sequence (\ref{415}) because
$\depth_{I(Z)} H_1 \ge 2$ implies $\Hom (K_D(n+c+1),H_1)\cong
H^0_{*}(U, {\tilde K_D(n+c+1)}^*\otimes {\tilde H_1})$ and the
right most map in the exact sequence $$ 0\rightarrow  \wedge ^3
(F(-\ell _2))\rightarrow \wedge ^3 (G_2(-\ell _2))\rightarrow
 \wedge ^2
(F(-\ell _2))\rightarrow H_2\rightarrow 0$$ (see \cite{AH}) must
correspond to the map $ \wedge ^2 (F(-\ell _2))\otimes D
\rightarrow H_2$ in (\ref{415}) and the later is surjective (which
one may prove directly as well, by applying $H^0_{*}(U,{\tilde
K_D}^*\otimes  (.))$ to (\ref{418}), to see $H^0_{*}(U,{\tilde
K_D}^*\otimes  {\tilde H_1})=0$).

To see that (\ref{416}) is exact we dualize (\ref{417}) and we get
$$0\rightarrow (I_D/I^2_D)^*\rightarrow G_2(-\ell _2)^*\otimes
D\rightarrow H_1^*\rightarrow 0$$ because $\Ext^1_D(
I_D/I^2_D,D)\cong \Ext^1_D( (I_D/I^2_D)\otimes K_D,K_D)=0$ by the
Cohen-Macaulayness of $(I_D/I^2_D)\otimes
K_D(n+c+1)\cong(I_D/I^2_D)^*$, cf. the proof of Theorem \ref{MAIN}
for the last isomorphims and \cite{KMMNP}, Ch. 6, for the
Cohen-Macaulayness. Applying $\Hom_D(.,K_D^*)$ to the last exact
sequence we get (\ref{416}) because $\depth _{I(Z)}D\ge 3$ implies
$\Hom _D(H_1^*,K_D^*)\cong \Hom _D(K_D,H_1)$ and
$\Ext^1_D(H_1^*,K_D^*)\cong H^1_{*}(U, \cH
om(\tilde{K_D},\tilde{H_1}))=0$ where the vanishing is due to the
Cohen-Macaulayness of $\Hom(K_D,H_1)$, which holds because we
already have proved the exactness of (\ref{415}). This concludes
the proof.
\end{proof}

\begin{corollary} \label{cod4}
 Let
$W(\underline{b};\underline{a})$ be the locus of good
determinantal schemes in $\PP^{n+c}$ where $n\ge 1$ and $c= 4$. If
$a_{i-min(3,t)}\ge b_{i}$ for $min(3,t)\le i \le t$, then $$ \dim
W(\underline{b};\underline{a})= \lambda _4+K_3+K_4.$$
\end{corollary}
\begin{proof} Due to Remark \ref{dep} and the assumption
$a_{i-min(3,t)}\ge b_{i}$ for $min(3,t)\le i \le t$, the set
$Z_2=Sing(X_2)$ satisfies $\depth _{I(Z_2)}D_2\ge 4$ by choosing
$C=Proj(A)$ general enough in $W(\underline{b};\underline{a})$.
Hence combining (\ref{new?}), Lemma \ref{key4} and Theorem
\ref{MAIN}, we are done.
\end{proof}

To apply Theorem \ref{MAIN} in the codimension $c=5$ case, it
suffices to prove that $$\Ext^1_{D_3}(I_{D_3}\otimes
I_3^*,I_3)_{\nu }= \Ext^1_{D_3}(I_{D_3}\otimes
M_3,M_3^*(-a_{t+2}-a_{t+3}))_\nu=0$$ for $\nu \le 0$. Since
$\depth_{I(Z_3)}D_3\ge 3$ and $I_3$ is a maximal Cohen-Macaulay
$D_3$-module, we have by (\ref{NM})

$$ \Ext^1_{D_3}(I_{D_3}\otimes M_3,M_3^*)\cong H^1_{*}(U_3,\cH om
(\cI _{X_3}\otimes \tilde{M}_3,\tilde{M}_3^*)) \cong $$
$$H^1_{*}(U_3,\cH om (\cI _{X_3}\otimes
S_2(\tilde{M_3)},\cO_{X_3}))\cong \Ext^1_{D_3}(I_{D_3}\otimes
S_2(M_3),D_3)$$ where $U_3=X_3-Z_3$. Since by Proposition 2.2(iii)
$K_{D_3}(n+c+1-\ell _3)\cong S_2(M_3)$ with $\ell _3 =\sum
_{j=0}^{t+1}a_j-\sum _{i=1}^tb_i$, we get (letting $B:=D_3$)
\begin{equation} \label{419} \Ext^1_{D_3}(I_{D_3}\otimes
I_3^*,I_3)=\Ext^1_B(I_B\otimes K_B(n+1+c-\ell
_3),B(-a_{t+2}-a_{t+3})) \end{equation} $$\cong
\Ext^1_B(I_B\otimes K_B(n+1+c),B)(\ell _3-a_{t+2}-a_{t+3}).$$

\begin{lemma} \label{411} Let $R\rightarrow B=R/I_B$ be a
codimension 3 good determinantal quotient, let $X \hookrightarrow
\PP^{n+c}$ be the corresponding embedding, and let $Z=Sing(X)$.
\begin{itemize}
\item[(a)] If $\depth _{I(Z)}B\ge 4$  then there is
an exact sequence $$ 0\rightarrow \Ext ^1_B(I_B\otimes
K_B(n+c+1),B)\rightarrow I_B/I_B^2\rightarrow (I_B/I_B^2)^{**}$$
which preserves the grading. In particular,
\begin{itemize}
\item[(a1)] $\Ext ^1_B(I_B\otimes
K_B(n+c+1),B)(\ell _3-a_{t+2}-a_{t+3})_{\nu }=0$ for
$\nu<a_{t+3}+a_{t+2}-a_{t+1}-a_t.$
\item[(a2)] If $Char(k)=0$, then $\Ext ^1_B(I_B\otimes
K_B(n+c+1),B)(\ell _3-a_{t+2}-a_{t+3})_{\nu }=0$ for $\nu\le
a_{t+3}+a_{t+2}-a_{t+1}-a_t.$
\end{itemize}
\item[(b)]If $\depth _{I(Z)}B\ge 5$ then there is an exact sequence $$
I_B/I_B^2\rightarrow (I_B/I_B^2)^{**}\rightarrow \Ext
^2_B(I_B\otimes K_B(n+c+1),B)\cong
H^1_{I(Z)}(I_B/I_B^2)\rightarrow 0$$ which preserves the grading.
\end{itemize}
\end{lemma}
\begin{remark} \label{rem412} Note that (a2) shows the desired vanishing because
 we in Theorem \ref{MAIN} have assumed $a_0\le a_1 \le \cdots \le
 a_{t+3}$.
\end{remark}
\begin{proof} (a) The Eagon-Northcott resolution associated to
$\varphi _3:F\rightarrow G_3=\oplus _{j=0}^{t+1} R(a_j)$ leads to
\begin{equation}\label{420} 0\rightarrow F_3:=\wedge
^{t+2}G_3^*\otimes S_2F\otimes \wedge ^tF\rightarrow F_2:=\wedge
^{t+1}G_3^*\otimes S_1F\otimes \wedge ^tF\rightarrow
\end{equation}

$$\rightarrow F_1:=\wedge ^{t}G_3^*\otimes \wedge ^tF\rightarrow
I_B\rightarrow 0.$$ Applying $\Hom_R(.,R)$ we get the exact
sequence $$ 0\rightarrow R\rightarrow F_1^*\rightarrow
F_2^*\rightarrow F_3^*\rightarrow \Ext^2_R(I_B,R)\cong
K_B(n+1+c)\rightarrow 0. $$ The tensorialization with $.\otimes
_RB$ leads to a complex
\begin{equation}\label{421}
0 \rightarrow(I_B/I_B^2)^* \rightarrow F_1^*\otimes B\rightarrow
F_2^*\otimes B\stackrel {\psi }{ \longrightarrow} F_3^*\otimes B
\rightarrow K_B(n+c+1)\rightarrow 0
\end{equation}
which is exact except in the middle where we have the homology
$I_B \otimes K_B(n+c+1)\cong \Tor _1^R(K_B(n+c+1),B)$. Indeed this
easily follows from the right exactness of $.\otimes _B B$ and the
left exactness of $\Hom_R(.,B)$ (applied to (\ref{420})). Since
(\ref{420}) also implies
\begin{equation}\label{422} 0\rightarrow H'_1:=\ker (\rho )\rightarrow
F_1\otimes _R B \stackrel {\rho }{ \longrightarrow}
I_B/I_B^2\rightarrow 0
\end{equation}
(observe that ${H'}_1$ is quite close to the 1. Koszul homology
$H_1$). By \cite{KP1}, Lemma 35, we have
$\depth_{\goth{m}}(I_B/I_B^2)^*\ge \depth _{\goth{m}}B-1$ and
hence by (\ref{NM}), $$\Ext^1_B(I_B/I_B^2,B)=0.$$

Dualizing (\ref{422}), it follows that
\begin{equation}\label{423}
0\rightarrow (I_B/I_B^2)^*\rightarrow F_1^*\otimes B\rightarrow
{H'}_1^*\rightarrow 0
\end{equation}
(and, if desirable, one may see $H_1^*\cong {H'}_1^*$). Since we
know the homology "in the middle" of (\ref{421}), we get the exact
sequences

\begin{equation}\label{424}
0\rightarrow {H'}_1^*\rightarrow \ker (\psi)\rightarrow I_B\otimes
K_{B}(n+c+1)\rightarrow 0,
\end{equation}
\begin{equation}\label{425} 0\rightarrow \ker (\psi) \rightarrow
F_2^*\otimes B\stackrel {\psi }{ \longrightarrow} F_3^*\otimes B
\rightarrow K_B(n+c+1)\rightarrow 0.
\end{equation}

Now we have the set-up to prove that $\Ext ^1_B(I_B\otimes
K_B(n+c+1),B)\cong \cK :=\ker (I_B/I_B^2\rightarrow
(I_B/I_B^2)^{**})$. Firstly note that dualizing (\ref{423}) once
more and comparing with (\ref{422}) in obvious way, we see that
$$0\rightarrow {H'}_1\rightarrow {H'}_1^{**}\rightarrow \cK
\rightarrow 0$$ by the snake-lemma. Now we apply $\Hom (.,B)$ to
(\ref{424}) and the left part of (\ref{425}). We get a commutative
diagram
\[
\begin{array}{cccccccc} F_2\otimes B &  \rightarrow & {H'}_1 &
\rightarrow & 0
\\   \downarrow  &  & \downarrow  \\
\Hom( \ker(\psi),B) & \rightarrow &
 \Hom({H'}^{*}_1,B) &  \rightarrow &   \Ext ^1_B(I_B\otimes K_B(n+c+1),B) & \rightarrow
 \Ext ^1_B (\ker (\psi),B)
\\    \downarrow  & &
 &     \\    \Ext ^1_B(\im (\psi),B) &  &
\end{array}
\]
from which we deduce the exact sequence
\begin{equation}\label{426}
 \Ext ^1_B (\im (\psi),B)\rightarrow \cK \rightarrow  \Ext ^1_B(I_B\otimes K_B(n+c+1),B) \rightarrow
 \Ext ^1_B (\ker (\psi),B).
\end{equation}

Hence it suffices to show $$\Ext ^1_B (\im (\psi),B)=0=\Ext ^1_B
(\ker (\psi),B).$$ By (\ref{425}) we have $$\Ext ^1_B (\im
(\psi),B)(n+c+1)\cong \Ext ^2_B (K_B,B) \cong \Ext ^2_B(K_B\otimes
K_B,K_B)$$
 $$\Ext ^1_B (\ker (\psi),B)(n+c+1)\cong \Ext ^3_B (K_B,B) \cong \Ext
^3_B(K_B\otimes K_B,K_B)$$ where the rightmost isomorphism is a
consequences of the spectral sequence used in \cite{her}, Satz 1.2
because we have $\depth_{I(Z)} B\ge 3$. By \cite{cm}; Corollary
3.4, we know that $\depth_{\goth{m}}S_2(K_B)\ge \depth_{\goth{m}}
B-1$. Hence by Gorenstein duality $\Ext ^{i}_B(S_2(K_B),K_B)=0$
for $i\ge 2$. Defining $\wedge $ by $$ 0\rightarrow \wedge
\rightarrow K_B\otimes K_B\rightarrow S_2(K_B)\rightarrow 0$$ and
noting that $\tilde{\wedge }|_{Proj(B)-Z}=0$, we get $\Ext ^{i}
_B(\wedge, K_B)=0$ for $i\le 3$ by (\ref{NM}) and the assumption
$\depth _{I(Z)}B \ge 4$. Combining we get $\Ext ^{i}_B(K_B\otimes
K_B,K_B)\cong \Ext ^{i}_B(S_2(K_B),K_B)=0$ for $i=2$ and $3$ as
required, i.e. $\cK\cong \Ext ^1_B(I_B\otimes K_B(n+c+1),B)$ by
(\ref{426}).

Now it is a triviality to see (a1) because the smallest  degree of
a minimal generator of $I_B$ is $\ell _3 -a_t-a_{t+1}.$

\vskip 2mm
 (a2) Since $\depth _{I(Z)}B\ge 2$,
we get $(I_B/I_B^2)^{**}\cong H^0(X-Z,\cI _X/\cI _X^2)$ and hence
that $\cK $ is isomorphic to $H^0_{I(Z)}(B)$. Similarly we prove
that the kernel of the "universal" derivation
$d:I_B/I_B^2\rightarrow\Omega_{R/k}\otimes _R B$ is
$H^0_{I(Z)}(B)$ which by \cite{EM}, Theorem 3, is isomorphic to
$I_B^{(2)}/I_B^2$ where $I_B^{(2)}$ is the second symbolic power
of $I_B$. Hence we have a grading-preserving isomorphism
\begin{equation}\label{427}
\Ext^1 _B(I_B\otimes K_B(n+c+1),B)\cong I_B^{(2)}/I_B^2.
\end{equation}

Now, in characteristic zero, $I_B^{(2)}\subset \goth{m} I_B$ by
\cite{EM}; Proposition 13, which shows that the smallest degree of
the minimal generators of $I_B^{(2)}$ is at least one less than
the smallest degree of the generators of $I_B$, i.e. we have
$$(I_B^{(2)})_{\ell _3-a_{t+2}-a_{t+3}+\nu}=0 \mbox{ for } \nu \le
0$$ and we conclude by (\ref{427}).

\vskip 2mm (b) Again since  $\depth _{I(Z)}B\ge 2$, we have
$(I_B/I^2_B)^{**}\cong H^0_*(X-Z,\cI _X/\cI ^2_X)$ and hence
$\coker[I_B/I^2_B\rightarrow (I_B/I^2_B)^{**}]\cong
H^1_{I(Z)}(I_B/I^2_B)\cong H^2_{I(Z)}(H'_1)$, cf. (\ref{422}) for
the last isomorphism. Using (\ref{424}), we get the exact sequence
$$ \Ext^1 _B(\ker(\psi ),B)\rightarrow \Ext
^1_B(H^{'*}_1,B)\rightarrow \Ext ^2_B(I_B\otimes
K_B(n+c+1),B)\rightarrow \Ext ^2_B(\ker(\psi ),B).$$ As argued in
(\ref{426}) and after (\ref{426}), we see that $\Ext
^1_B(\ker(\psi ),B)=0$ and $$\Ext ^2_B(\ker(\psi ),B)(n+c+1)\cong
\Ext ^4_B(K_B\otimes K_B,K_B)\cong \Ext ^4_B(S_2(K_B),K_B)=0$$
where the last isomorphism to the second symmetric power follows
from the fact that $\depth _{I(Z)}B\ge 5$ implies $\Ext
^{i}_B(\wedge ,K_B)=0$ for $i\le 4$, and the vanishing to the
right follows from $\depth _{\goth{m}} S_2(K_B)\ge
\depth_{\goth{m}}B-1$. Since by (\ref{NM}), $$\Ext
^1_B(H^{'*}_1,B)\cong H^1_*(U,\cH
om(\tilde{H^{'*}_1},\tilde{B}))\cong H^1_*(U,\tilde{H'_1})\cong
H^2_{I(Z)}(H'_1)$$ we are done.
\end{proof}

 \begin{remark} For generic determinantal schemes one knows that
 $\depth _{I(Z)}(I_B/I^2_B)\ge 2$ by \cite{b-v}. So the vanishing
 of $\Ext ^1_B(I_B\otimes K_B(n+c+1),B)_{\nu }$ under reasonable
 genericity assumptions is expected (for any $\nu $).
 \end{remark}

 \begin{corollary} \label{cod5} Let
$W(\underline{b};\underline{a})$ be the locus of good
determinantal schemes in $\PP^{n+c}$ where $n\ge 1$ and $c= 5$. If
$a_{i-min(3,t)}\ge b_{i}$ for $min(3,t)\le i \le t$  and
$Char(k)=0$, then $$ \dim W(\underline{b};\underline{a})= \lambda
_5+K_3+K_4+K_5.$$
 \end{corollary}
\begin{proof} It follows from Remark \ref{dep} and the assumption
$a_{i-min(3,t)}\ge b_{i}$ for $min(3,t)\le i \le t$ that the set
$Z_j=Sing(X_j)$  has $\depth _{I(Z_j)}D_j \ge 4$ for $j=2$ and 3
provided $C=Proj(A)$ is chosen general enough in
$W(\underline{b};\underline{a})$. By (\ref{new?}), Lemma
\ref{key4}, (\ref{419}), Lemma \ref{411}(a2) and Remark
\ref{rem412}
 the assumptions of Theorem
\ref{MAIN} are fulfilled  and we conclude by applying it.
\end{proof}

Now we state the last Corollaries of this section which shows that
the upper bound of $\dim W(\underline{b};\underline{a})$ given in
Theorem \ref{upperbound} is indeed equal to $\dim
W(\underline{b};\underline{a})$ for {\em all} $c\ge 3$ and most
values of $a_0,a_1,\cdots ,a_{t+c-2};b_1, \cdots ,b_t$. Our result
is based upon Remark \ref{rem44} and the proof of Theorem
\ref{MAIN}. Indeed, we have seen that $$ \dim
W(\underline{b};\underline{a})= \lambda _c+K_3+K_4+ \cdots +K_c$$
provided $$\dim W(F,G_{c-1})= \lambda_{c-1}+K_3+K_4+\cdots
+K_{c-1}$$ \noindent and
\begin{equation}
\label{428} \Hom _R(I_{D_{c-2}},I_{c-1})_0=0.
\end{equation}

 \begin{corollary} \label{cod6}
 Let
$W(\underline{b};\underline{a})$ be the locus of good
determinantal schemes in $\PP^{n+c}$ where $n\ge 0$ and $c\ge 6$.
 Assume $a_{i-min(3,t)}\ge b_{i}$ for $min(3,t)\le i \le t$, $Char(k)=0$ and

$(i_6): \quad a_{t+4}>a_{t-1}+a_t+a_{t+1}+a_{t+2}-a_0-a_1-a_2$,

$(i_7): \quad
a_{t+5}>a_{t-1}+a_t+a_{t+1}+a_{t+2}+a_{t+3}-a_0-a_1-a_2-a_3$,

....

$(i_c): \quad a_{t+c-2}>\sum _{j=t-1}^{t+c-4}a_{j}-\sum
_{j=0}^{c-4}a_j$.

 Then, $$ \dim
W(\underline{b};\underline{a})= \lambda _c+K_3+\cdots +K_c.$$
 \end{corollary}
\begin{proof} By the Eagon-Northcott resolution the largest
possible degree of a generator of $I_{D_{c-2}}$ is $\ell
_c-\sum_{j=0}^{c-4}a_j-a_{t+c-3}-a_{t+c-2}$ where $\ell
_c=\sum_{j=0}^{t+c-2}a_j-\sum_{i=0}^tb_i$ and the smallest
possible degree of a generator of $I_{c-1}\cong
I_{D_c}/I_{D_{c-1}}$ is $\ell _c-\sum _{j=t-1}^{t+c-3} a_{j}$
because $a_0\le a_1\le \cdots \le a_{t+c-2}$. Hence if the latter
is strictly larger than the former, i.e. if $$a_{t+c-2}>\sum
_{j=t-1}^{t+c-4}a_{j}-\sum _{j=0}^{c-4}a_j$$ then
$\Hom(I_{D_{c-2}},I_{c-1})_0=0$ and we conclude using the argument
of (\ref{428}) and Corollaries \ref{cod3}, \ref{cod4} and
\ref{cod5}.
\end{proof}

\begin{remark}\label{debil}
\label{weak}(1) If we want to skip the characteristic zero
assumption, we can avoid the use of Corollary \ref{cod5} by
introducing the assumption $$(i_5): \quad
a_{t+3}>a_{t-1}+a_t+a_{t+1}-a_0-a_1.$$ We still get
 $ \dim
W(\underline{b};\underline{a})= \lambda _c+K_3+\cdots +K_c$,
supposing $a_{i-min(3,t)}\ge b_{i}$ for $min(3,t)\le i \le t$ and
$(i_5)$, $(i_6)$,...,$(i_c)$.

(2) We can further weaken $a_{i-min(3,t)}\ge b_{i}$ for
$min(3,t)\le i \le t$ to $a_{i-min(2,t)}\ge b_{i}$ for
$min(2,t)\le i \le t$ by avoiding Corollary \ref{cod4} and
assuming in addition $(i_4): \quad a_{t+2}>a_{t-1}+a_t-a_0.$
 \end{remark}

 \begin{remark} \label{dim0}
 While Corollaries \ref{cod3}, \ref{cod4} and
\ref{cod5} do not apply to the case when $
W(\underline{b};\underline{a})$ is the locus of zero-dimensional
determinantal schemes, Corollary \ref{cod6} and Remark \ref{weak}
do apply to the zero-dimensional case. In particular, using Remark
\ref{weak} (1) (resp. (2)) for $c=5$ (resp. $c=4$), we get a
single assumption, namely $(i_5)$ (resp. $(i_4)$) in addition to
$a_{i-min(3,t)}\ge b_{i}$ for $min(3,t)\le i \le t$ (resp.
$a_{i-min(2,t)}\ge b_{i}$ for $min(2,t)\le i \le t$) which
suffices for having $\dim W(\underline{b};\underline{a})$ equal to
the upper bound given in Theorem \ref{upperbound} for the zero
schemes as well.
 \end{remark}

It is worthwhile to point out that this last remark on
zero-schemes works also in the codimension $c=3$  case, and here
the $(i_3)$ assumption is very weak. We have

\begin{corollary} \label{cod3dim0} Let
$W(\underline{b};\underline{a})$ be the locus of good
determinantal schemes  in $\PP^{n+3}$ of codimension 3. If
$a_{i-min(2,t)}\ge b_{i}$ for $min(2,t)\le i \le t$ and if in
addition $$(i_3): \quad a_{t+1}>a_{t-1} $$ then $\dim
W(\underline{b};\underline{a})=\lambda _3+K_3.$\end{corollary}
\begin{proof} Slightly extending  Remark \ref{dep} by introducing the determinantal
hypersurface $X_1=Proj(D_1)$ we have $\depth _{I(Z_1)}D_1 \ge 3$
and $\depth _{I(Z_2)}D_2\ge 2$ by choosing $C=Proj(A)\in
W(\underline{b};\underline{a})$ general enough. It follows that
(\ref{47}) is exact also for $i=1$, and since $a_{i-min(2,t)}\ge
b_{i}$ for $min(2,t)\le i \le t$ implies $\Hom
_R(I_{D_1},I_2)_0=0$ we get $$\hom (I_{D_2},I_2)_0\cong \dim
(D_2)_{a_t-a_{t+1}}.$$

Hence Proposition \ref{main0} (ii) for $c=3$ applies to explicitly
get a lower bound of $\dim W(\underline{b};\underline{a})$, which
combining (\ref{Mc0}) and (\ref{lambdas}) turns out to be $\lambda
_3+K_3.$ Hence, $\dim W(\underline{b};\underline{a})=\lambda
_3+K_3$ by Theorem \ref{upperbound}.
\end{proof}

\section{Unobstructedness of determinantal schemes}

In this section we keep the notation introduced in sections 3 and
4 and we consider the problem of when the closure of $
W(\underline{b};\underline{a})$ is an irreducible component of
$\Hi ^p(\PP^{n+c})$ and when $\Hi ^p(\PP^{n+c})$ is smooth or, at
least, generically smooth along $ W(\underline{b};\underline{a})$.

\vskip 4mm Through this section we {\em always} assume $n\ge 1$
and $c\ge 2$, The following result is crucial to our work in this
section:

\begin{theorem}\label{main5} Let $C\subset \PP^{n+c}$ be a good
determinantal scheme of dimension $n\ge 1$, let $C=X_c\subset
X_{c-1}\subset ...\subset X_2\subset \PP^{n+c}$ be the flag
obtained by successively deleting columns from the right hand side
and let $Z_i\subset X_i$ be some closed subset such that
$X_i-Z_i\subset \PP^{n+c}$ is a local complete intersection.

(i) If $\depth _{I(Z_{i})}D_i \ge 3 $ for $2\le i \le c-1$, and if
$\Ext ^1_{D_i}(I_{D_i}/I^2_{D_i},I_i)_0\hookrightarrow \Ext
^1_{D_i}(I_{D_i}/I^2_{D_i},D_i)_0$ for $i=2,...,c-1,$ then $C$
(and each $X_i$) is unobstructed, and $$\dim_{X_{i+1}}\Hi
^{p_{i+1}}(\PP^{n+c})=\dim_{X_{i}}\Hi ^{p_{i}}(\PP^{n+c})+\dim
(N_{D_{i+1}/D_{i}})_0-\hom (I_{D_{i}},I_{i})_0$$ \noindent for
$i=2,3,...,c-1$.

 (ii) If $a_0\le a_1 \le ... \le a_{t+c-2}$, $b_1\le ... \le b_t$
 and $a_{i-min(2,t)}\ge b_{i}$ for $min(2,t)\le i \le t$, and
 if  a sufficiently general $C\in W(\underline{b};\underline{a})$
  satisfies  $$\Ext ^1_{D_i}(I_{D_i}/I^2_{D_i},I_{i})_0=0 \mbox{ for }
 i=2,...,c-1$$
 then  $\overline{
W(\underline{b};\underline{a})}$ is a generically smooth
irreducible component of $\Hi ^p(\PP^{n+c})$.
\end{theorem}

\begin{proof} (i) First of all we claim that there are 2 short exact
sequences, the vertical and the horizontal one, fitting into a
commutative diagram (whose square is cartesian)

\begin{equation}\label{51bis}
\begin{array}{ccccccccc} & & &  & & 0 \\
& & & & &   \downarrow \\ &  & & & & \Hom_R(I_{D_i},I_i)_0 \\ & &
& & &   \downarrow \\ & & & \hskip 10mm A^1 & \stackrel {T_{pr_2}
}{ \longrightarrow} & \Hom_R(I_{D_i},D_i)_0
\\ & &  &  T_{pr_1}\downarrow  &  \smallbox  & \downarrow \\ 0 \rightarrow & \Hom
(I_i,D_{i+1})_0 & \rightarrow & \Hom (I_{D_{i+1}},D_{i+1})_0 &
\rightarrow & \Hom_R(I_{D_i},D_{i+1})_0 &\rightarrow &  0 \\ & & &
& & \downarrow  \\
 & & &  & & 0
\end{array}
\end{equation}
where $A^1$ is the tangent space of the Hilbert flag scheme
$D(p_{i+1},p_i)$ at $(X_{i+1}\subset X_i)$ and $T_{pr_i}$ the
tangent maps of the projections $pr_i$ (see Proposition 4.1 for
details). Since the vertical sequence is exact by assumption and
the tangent space description of the Hilbert flag scheme and its
projections are well known (\cite{KMMNP}, Chapter 6, and note that
the zero piece of the graded $\Hom $'s above and the corresponding
global sections of their sheaves of \cite{KMMNP} coincide by
(\ref{NM})), we only have to prove the short-exactness of the
horizontal sequence. Hence it suffices to prove that $T_{pr_2}$ is
surjective. To see it, it suffices to slightly generalize the
argument in the proof of Proposition 4.1 where we showed that the
dimension of the fiber is $\ge m_c(0)$. We skip the details since
\cite{KMMNP}, Theorem 10.13 shows more. Indeed, it contains a
deformation theoretic argument which shows that $pr_2$ is not only
dominating but also "infinitesimal dominating or  surjective"
(i.e. smooth at $(X_{i+1}\subset X_i)$). In particular, we have
that the tangent map $T_{pr_2}$ is surjective, cf. Remark
\ref{52NEW} for another argument.

  By the
proof of Theorem 10.13 of \cite{KMMNP} (see Remark \ref{52NEW} for
an easy argument), $D(p_{i+1},p_i)$ is smooth at $(X_{i+1}\subset
X_{i})$ provide $\Hi ^{p_{i}}(\PP^{n+c})$ is smooth at $X_i$.
Since the tangent map of the first projection $pr_1:
D(p_{i+1},p_i)\rightarrow \Hi ^{p_{i+1}}(\PP^{n+c})$ is
surjective, we get that $\Hi ^{p_{i+1}}(\PP^{n+c})$ is smooth at
$X_{i+1}$. By induction $C$ (and each $X_i$) is unobstructed since
$X_2$ is unobstructed \cite{elli}, and the two exact sequences of
(\ref{51bis}) easily lead to the  dimension of $\dim _{X_{i+1}}
\Hi ^{p_{i+1}} (\PP^{n+c})$  because $N_{D_{i+1}/D_{i}}=\Hom
(I_i,D_{i+1})$.

(ii) To prove that $\overline{W(\underline{b};\underline{a})}$ is
an irreducible component, we use the notation of Proposition 4.1
and we may by induction suppose  that $\overline{W(F,G_{c-1})}$ is
an irreducible component of $\Hi ^{p_{c-1}}(\PP^{n+c})$ since
$\overline{W(F,G_2)}$ is an irreducible component by \cite{elli}.
We have
\begin{equation}\label{51} \dim D(F,G_c,G_{c-1})\ge \dim
W(F,G_{c-1})+m_c(0)
\end{equation}
by Proposition 4.1 (ii) while for an irreducible component $V$ of
$D(p_c,p_{c-1})$ containing $D(F,G_c,G_{c-1})$ we must have
\begin{equation}\label{52}
\dim  V\le \dim W(F,G_{c-1})+\dim(N_{D_{c}/D_{c-1}})_0
\end{equation}
because $\dim(N_{D_{c}/D_{c-1}})_0$ is the fiber dimension of
$pr_2$ at $(X_c\subset X_{c-1})$. Since $\depth
_{I(Z_{c-1})}D_{c-1}\ge 3$, we have by (\ref{DiMi}) $$\dim
(N_{D_c/D_{c-1}})_0=m_c(0).$$

Combining (\ref{51}) and (\ref{52}) we get $\dim
D(F,G_c,G_{c-1})\ge \dim V$ and hence
$\overline{D(F,G_c,G_{c-1})}$ is
 an irreducible component of $D(p_c,p_{c-1})$. Since
 the first projection $pr_1:D(p_c,p_{c-1})\rightarrow\Hi
^{p_c}(\PP^{n+c})$ is smooth at $(X_{i+1}\subset X_i)$ by the
surjectivity of $T_{pr_1}$ and the smoothness of $D(p_{i+1},p_i)$
at $(X_{i+1}\subset X_i)$, we get that
 $\overline{W(\underline{b};\underline{a})}$ is an
irreducible component, which necessarily is generically smooth
because  $\Hi ^{p_c}(\PP^{n+c})$ is smooth at a general point $C$
by the first part of the proof  and by Remark \ref{dep}.
\end{proof}

\begin{remark}\label{52NEW}
{\rm If, in Theorem \ref{main5}, we suppose $\depth
_{I(Z_{i})}D_i\ge 4$  we may easily see the surjectivity of
$T_{pr_2}$ in the following way. Using (\ref{NM}), we get $\Ext
^1_{D_{i}}(I_i,R_i)=\Ext ^2_{D_{i}}(I_i,I_i)=0$ by the depth
condition above. Applying  $\Hom _{D_{i}}(I_i,.)$ to  the exact
sequence $0\rightarrow I_i\rightarrow D_i\rightarrow D_{i+1}
\rightarrow 0$, we get $\Ext ^1_{D_{i}}(I_i,D_{i+1})=0$ and the
lower horizontal sequence of (\ref{51bis}) is short exact and we
easily conclude. Finally using the vanishing of $\Ext
^1_{D_{i}}(I_i,D_{i+1})_0$, it follows from (\ref{NM}) that
$H^1(U_i,\tilde{N}_{D_{i+1}/D_{i}})\cong \Ext
^1_{D_{i+1}}(I_i/I^2_i,D_{i+1})_0=0$. Then it is not difficult to
see that $D(p_{i+1},p_i)$ is smooth at $(X_{i+1}\subset X_i)$
provided $\Hi ^{p_{i}}(\PP^{n+c})$ is smooth at $X_i$.}

\end{remark}

\vskip 2mm To apply  Theorem \ref{main5}(ii) in the codimension
$c=3$ case, it suffices to prove that $$\Ext
^1_{D_2}(I_{D_2}/I_{D_2}^2,I_2)_0=0.$$ By (\ref{NM}) and
(\ref{412}) we see that ($U_2=X_2-Z_2$)
\begin{equation} \label{53}\Ext
^1_{D_2}(I_{D_2}/I_{D_2}^2,I_2)\cong H^1_*(U_2,\cH
om(\cI_{X_2},\tilde{K}_{D_2}(n+4),\cO_{X_2})(\ell _2-a_{t+1}))
\end{equation}
$$\cong H^1_*(U_2,\cI_{X_2}/\cI_{X_2}^2(\ell _2-a_{t+1}))$$ and we
consider two cases:

If $\depth_{I(Z_2)}D_2\ge 4$, we get $\depth
_{I(Z_2)}I_{D_2}/I_{D_2}^2\ge 3$ \cite{AH} and the group in
(\ref{53}) vanishes.

If $\depth_{I(Z_2)}D_2=3$ (e.g. $X_2$ is smooth and 2
dimensional), the group, in degree zero, is clearly $
H^1(U_2,\cI_{X_2}/\cI_{X_2}^2(\ell _2-a_{t+1}))$, and we have to
suppose it vanishes in order to conclude that
$\overline{W(\underline{b};\underline{a})}$ is a generically
smooth component of $\Hi ^p(\PP^{n+3})$ of dimension $\lambda
_3+K_3$. All this is essentially \cite{KMMNP}, Corollary 10.15
(ii).

\vskip 3mm The case $c=4$ is also straightforward. In this case it
suffices to see that (\ref{53}) vanishes and that
\begin{equation} \label{54}\Ext
^1_{D_3}(I_{D_3}/I_{D_3}^2,I_3)=0.
\end{equation}

If we suppose
\begin{equation}\label{55} \depth _{I(Z_2)} D_2 \ge 4 \mbox{ and } \depth _{I(Z_3)} D_3 \ge 4
\end{equation}
we claim that both group vanish. We only need to prove (\ref{54}).
Since $\depth _{I(Z_2)} D_2 \ge 4$ it follows from (\ref{NM}) that
(\ref{47}) is short-exact for $i=3$. Using Lemma \ref{key4} and
(\ref{new?}), we see that (\ref{48}) is short-exact for $i=3$ as
well, i.e. we have exact sequences

\begin{equation}\label{56} \hspace{10mm}
0\rightarrow D_3(a)\rightarrow \Hom _R(I_{D_3},I_3)\rightarrow
\Hom _R(I_{D_2},I_3)\rightarrow 0
\end{equation}
$$ \hspace{60mm} \|$$ $$0\rightarrow \Hom(I_{D_2}\otimes
I_2^*,I_2(a))\rightarrow \Hom_R (I_{D_2},I_2(a))\rightarrow
\Hom_{D_2}(I_{D_2}\otimes I^*_2,D_3(a))\rightarrow 0 $$ where
$a=a_{t+1}-a_{t+2}$. By Lemma \ref{key4}, the codepth of $\Hom
(I_{D_2}\otimes I^*_2,I_2(a))$ is at most 1 while (\ref{412})
shows the same conclusion for $\Hom (I_{D_2},I_2(a))$. The lower
exact sequence of (\ref{56}) therefore shows that the codepth of
$\Hom_{D_2} (I_{D_2}\otimes I^*_2,D_3(a))$ is at most 1 as a
$D_3$-module. The upper sequence shows that

\begin{equation}\label{57bis}
\depth _{\goth{m}} \Hom _R(I_{D_3},I_3)\ge \depth
_{\goth{m}}D_3-1. \end{equation}

\noindent Now since $\depth _{I(Z_3)}D_3\ge 4$, we get $\depth
_{I(Z_3)} \Hom_{D_3}(I_{D_3}/I^2_{D_3},I_3)\ge 3$ and hence by
(\ref{NM}) that (\ref{54}) holds. By Remark \ref{dep}, we see that
(\ref{55}) holds for a general $C\in
W(\underline{b};\underline{a})$ provided $a_{i-min(3,t)}\ge b_{i}$
for $min(3,t)\le i \le t$. Combining with Corollary \ref{cod3} and
Corollary \ref{cod4} we get

\vskip 2mm
\begin{corollary} \label{component34} Let
$W(\underline{b};\underline{a})$ be the locus of good
determinantal schemes in $\PP^{n+c}$ where $n\ge 2$ and $c=3$ or
$c=4$. If $a_{i-min(3,t)}\ge b_{i}$ for $min(3,t)\le i \le t$,
then $\overline{W(\underline{b};\underline{a})}$ is a generically
smooth, irreducible component of $\Hi ^p(\PP^{n+c})$ of dimension
$\lambda _c+K_3+...+K_c.$
\end{corollary}

\begin{remark} If $c=4$, then the assumption (\ref{55}) excludes
the interesting case when $W(\underline{b};\underline{a})$
parameterizes good determinantal curves in $\PP^5$. To consider
this case we will weaken (\ref{55}) and only suppose $\depth
_{I(Z_2)}D_2 \ge 4$. Recalling that (\ref{412}) leads to
$$H^1_*(U_2,\cH om(\cI _{X_2}\otimes \cI^*_2,\cO_{X_2}))\cong
H^1_*(U_2,\cI_{X_2}/\cI ^2_{X_2}(\ell _2-a_{t+1}))=0,$$ we have by
(\ref{47new}) injections

\begin{equation}\label{57}
\Ext ^1_{D_3}(I_{D_3}/I^2_{D_3},I_3)_0\hookrightarrow H^1(U_2,\cH
om(\cI_{X_{2}}\otimes \cI ^*_2,\cO_{X_3}(a))) \hookrightarrow
H^2(U_{2},\cH om(\cI_{X_{2}}\otimes \cI ^*_{2},\cI_{2}(a)))
\end{equation} $$ \hspace{25mm} \| \hspace{55mm} \|$$ $$
\hspace{35mm}
  H^1(U_{2},\cI_{X_{2}}/\cI^2_{X_2}\otimes \cO_{X_3}(\ell _2-a_{t+2}))
  \hookrightarrow
 H^2(U_{2},\cI_{X_{2}}/\cI^2_{X_2}\otimes \tilde{K} ^*_{D_2}(a'))
$$ where $a=a_{t+1}-a_{t+2}$ and $a'=2\ell _2-6-a_{t+1}-a_{t+2}$.
In the interesting case $C=X_4\subset X_3\subset X_2\subset \PP^5$
where $X_2$ is smooth, then $U_{2}=X_2$. In particular if one of
the groups of (\ref{57}) vanishes, then
$\overline{W(\underline{b};\underline{a})}$ is still  a
generically smooth, irreducible component of $\Hi ^p(\PP^{5})$ of
dimension $\lambda _4+K_3+K_4.$
\end{remark}

As a Corollary of the first part of Theorem \ref{main5} we also
get the unobstructedness and the vanishing of
$H^{i}_{*}(C,{\cN}_C)$ for good determinantal schemes $C\in
\PP^{n+c}$ of codimension $3\le c \le 4$. For $c=3$, the
unobstructedness is essentially proved in \cite{KMMNP}, Corollary
10.15 and the vanishing of $H^{i}_{*}(C,{\cN}_C)$ is shown in
\cite{KP1}, Lemma 35.

\begin{corollary}\label{new54} Let $C=Proj(A)\subset \PP^{n+c}$ be
a good determinantal scheme of dimension $n\ge 2$ for which there
is a flag satisfying $\depth  _{I(Z_{i})}(D_i)\ge 4$ for $2\le i
\le c-1.$ If $3\le c\le 4$, then $C$ is unobstructed and the
normal module $N_A:=Hom_R(I_A,A)$ satisfies $\depth _{\goth
{m}}(N_A)\ge n-1$. In particular $$ H^{i}_*(C,{\cN}_C)=0 \mbox{
for } 1\le i \le n-2.$$
\end{corollary}
\begin{proof} Due to the vanishing of (\ref{53}) and (\ref{54}),
the unobstructedness of $C$ follows at once from Theorem
\ref{main5}(i). Moreover exactly as we managed to show that the
exact sequences of (\ref{56}) implied (\ref{57bis}), we may see
that the graded exact horizontal and vertical sequences of
(\ref{51bis}) for $i=2$ (and not only the degree zero piece of
these sequences) imply

\begin{equation}\label{5??}
\depth _{\goth{m}} \Hom (I_{D_{3}},D_3)\ge  \dim D_3 -1
\end{equation}
\noindent because we have $\depth  _{I(Z_{2})}(D_2)\ge 4$ and
hence $\Ext ^1_{D_{2}}(I_{D_2}/I^2_{D_{2}},I_2)=0$ by (\ref{53}).
This shows what we want for $c=3$. Finally the same argument for
$i=3$, assuming (\ref{55}) and using (\ref{57bis}) and
(\ref{5??}), leads to $\depth _{\goth{m}} \Hom
_R(I_{D_{4}},D_4)\ge \dim D_4 -1$ and we are done.
\end{proof}

\vskip 2mm In the following example we will see that Corollary
\ref{component34} does not always work for determinantal curves
$C\subset \PP^5$, i.e.  the closure of
$W(\underline{b};\underline{a})$ is not necessarily an irreducible
component of $\Hi ^p(\PP^{5})$ although by Corollary \ref{cod4} we
know that $\dim W(\underline{b};\underline{a})$ is indeed $\lambda
_4+K_3+K_4$.

\begin{example} \label{ex1}
Let $C\subset \PP^5$ be a smooth good determinantal curve of
degree 15 and arithmetic genus 10 defined by the maximal minors of
a $3\times 6$ matrix with linear entries. The closure of
$W(\underline{b};\underline{a})$  inside $\Hi ^{15t-9}(\PP^5)$ is
not an irreducible  component. In fact, let $ H_{15,10}\subset \Hi
^{15t-9}(\PP^5)$ be the open subset parameterizing smooth
connected curves of degree $d=15$ and arithmetic genus $g=10$. It
is well known that any irreducible component of $H_{15,10}$ has
dimension $\ge 6d+2(1-g)=72$; while by Corollary \ref{cod4}, $\dim
W(\underline{b};\underline{a})=64.$
\end{example}

\vskip 2mm For the codimension $c=5$ case we have

\begin{corollary} \label{component5} Let
$W(\underline{b};\underline{a})$ be the locus of good
determinantal schemes in $\PP^{n+c}$ where $n\ge 1$ and $c=5$. If
$a_{i-min(3,t)}\ge b_{i}$ for $min(3,t)\le i \le t$  and if
$W(\underline{b};\underline{a})$ contains a determinantal scheme
$C=Proj(D_5)$ whose flag $R\rightarrow D_2\rightarrow
D_3\rightarrow D_4\rightarrow D_5$ obtained by deleting columns of
"largest possible degree" satisfies (with $Z_i=Sing(X_i)$),
$\depth _{I(Z_2)}D_2\ge 4$, $\depth _{I(Z_3)}D_3\ge 5$ and
$H^1(X_3-Z_3,\cI ^2_{X_3}(\ell _3-2a_{t+3}))=0$, then
$\overline{W(\underline{b};\underline{a})}$ is a generically
smooth, irreducible component of $\Hi ^p(\PP^{n+5})$ (of dimension
$\lambda _5+K_3+K_4+K_5$ provided $Char(k)=0$).
\end{corollary}
\begin{proof} First of all note that by Remark \ref{dep}, if we choose
$C\in W(\underline{b};\underline{a})$ general enough, we have
$\depth _{I(Z_2)}D_2 \ge 4$ and $\depth _{I(Z_3)}D_3 \ge 5$. By
Theorem \ref{main5} and the conclusion of (\ref{55}), it suffices
to show
\begin{equation}\label{59} \Ext
^1_{D_4}(I_{D_4}/I^2_{D_4},I_4)_0=0.
\end{equation}
 By Lemma \ref{ext}(b) and (\ref{54}) we must show that
 $\Ext ^2_{D_3}(I_{D_3}\otimes I_3^*,I_3)_{a_{t+2}-a_{t+3}}=0.$

Looking to (\ref{419}), this group is isomorphic to  $$\Ext
^2_{D_3}(I_{D_3}\otimes K_{D_3}(n+c+1),D_3)_{\ell _3-2 a_{t+3}}$$
which by Lemma \ref{411} (b) is further isomorphic to
$H^1_{I(Z_3)}(I_{D_3}/I^2_{D_3})_{\ell _3-2a_{t+3}}$. Since
$X_3=Proj(D_3)$ is Cohen-Macaulay, the cohomology sequence
associated to $$0\rightarrow I^2_{D_3}\rightarrow
I_{D_3}\rightarrow I_{D_3}/I^2_{D_3}\rightarrow 0$$ gives us
$$H^1_{I(Z_3)}(I_{D_3}/I^2_{D_3})_{\ell _3-2a_{t+3}} \cong
H^2_{I(Z_3)}(I^2_{D_3})_{\ell _3-2a_{t+3}}\cong H^1(X_3-Z_3,\cI
^2_{X_{3}}(\ell _3-2a_{t+3})) $$ and we get (\ref{59}).
\end{proof}

\vskip 2mm We will now give two examples. The first one will be a
smooth determinantal surface $S\subset \PP^7$ whose flag
$R\rightarrow D_2\rightarrow D_3\rightarrow D_4\rightarrow D_5$
obtained by deleting columns of "largest possible degree"
satisfies all hypothesis required in Corollary \ref{component5}
and hence $\overline{W(\underline{b};\underline{a})}$ is a
generically smooth, irreducible component of $\Hi ^p(\PP^{6})$ of
dimension $\lambda _5+K_3+K_4+K_5$ ($Char(k)=0$). The second one
will  be a smooth determinantal curve $C\subset \PP^6$ hence the
condition $\depth _{I(Z_3)}D_3\ge 5$ is not fulfilled and in this
case we will see that the closure of
 $W(\underline{b};\underline{a})$ is not  an
irreducible component of $\Hi ^p(\PP^{6})$ although by Corollary
\ref{cod5} we know that $\dim W(\underline{b};\underline{a})$ is
indeed $\lambda _4+K_3+K_4+K_5$ ($Char(k)=0$).

\begin{example}\label{ex2} (1) Let $S\subset \PP^7$ be a smooth
good determinantal surface of degree 6 defined by the maximal
minors of a $2\times 6$ matrix with general linear entries. Let
$R\rightarrow D_2\rightarrow D_3\rightarrow D_4\rightarrow D_5$ be
the flag  obtained by deleting  columns from the right hand side.
With the computer program  Macaulay \cite{Mac} we check that all
hypothesis required in Corollary \ref{component5} are satisfied,
i.e., $\depth _{I(Z_2)}D_2\ge 4$, $\depth _{I(Z_3)}D_3\ge 5$ and
$H^1(X_3-Z_3,\cI ^2_{X_3}(\ell _3-2a_{t+3}))=H^1(X_3,\cI
^2_{X_3}(2))=0$.
 By Corollary \ref{component5} the closure
of $W(\underline{b};\underline{a})$ inside $\Hi ^{p}(\PP^7)$ is a
generically smooth irreducible  component of dimension 57.

(2) Let $C\subset \PP^6$ be a smooth good determinantal curve of
degree 21 and arithmetic genus 15 defined by the maximal minors of
a $3\times 7$ matrix with linear entries. Since $\dim(C)=1$, we
have $\dim (X_3)=3$ and hence $\depth _{I(Z_3)}D_3\le 4$. The
closure of $W(\underline{b};\underline{a})$ inside $\Hi
^{21t-14}(\PP^6)$ is not an irreducible  component. In fact, let $
H_{21,15}\subset \Hi ^{21t-14}(\PP^6)$ be the open subset
parameterizing smooth connected curves of degree $d=21$ and
arithmetic genus $g=15$. It is well known that any irreducible
component of $H_{21,15}$ has dimension $\ge 7d+3(1-g)=105$; while
by Corollary \ref{cod5}, $\dim W(\underline{b};\underline{a})=90.$
\end{example}

Our final Corollaries are similar to Corollaries \ref{cod6} and
\ref{cod3dim0}. To apply the final part of Theorem \ref{main5} we
must show that $\Ext^1_{D_{i}}(I_{D_{i}}/I^2_{D_{i}},I_i)_0=0$ for
$i=2,...,c-1$. Using, however, the upper sequence of (\ref{48new})
it suffices to show that
$\Ext^1_{D_{i-1}}(I_{D_{i-1}}/I^2_{D_{i-1}},I_i)_0=0$ provided
$\depth _{I(Z_{i-1})}D_{i-1}\ge 4$. This vanishing is fulfilled if
\begin{equation}\label{510}
\Ext ^1_{R}(I_{D_{i-1}},I_i)_0=0 \mbox{ for } i=4,...,c-1
\end{equation}
(since $\Ext ^1_{D_i}(I_{D_i}/I^2_{D_i},I_i)_0=0$ for $i=2,3$
provided $\dim_{I(Z_{i})} D_i\ge 4$ by (\ref{53}) and (\ref{54})).

\begin{corollary}\label{component6}
Let $W(\underline{b};\underline{a})$ be the locus of good
determinantal schemes in $\PP^{n+c}$ where $n\ge 1$ and $c\ge 5$,
or $n\ge 2$ and $3\le c \le 4$. Assume $a_{i-min(3,t)}\ge b_{i}$
for $min(3,t)\le i \le t$
 and

$(j_5):$ $a_{t+3}>a_{t-1}+a_t+a_{t+1}-a_0-b_1,$

$(j_6):$ $a_{t+4}>a_{t-1}+a_t+a_{t+1}+a_{t+2}-a_0-a_1-b_1,$

....

$(j_c):$ $a_{t+c-2}>\sum _{j=t-1}^{t+c-4}a_{j}-\sum _{j=0}^{c-5}
a_j-b_1.$

Then $\overline{W(\underline{b};\underline{a})}$ is a generically
smooth, irreducible component of $\Hi ^p(\PP^{n+c})$ of dimension
$\lambda _c+K_3+...+K_c$.
\end{corollary}
\begin{proof} The relation of $I_{D_{c-2}}$ of the largest
possible degree is $\ell _c-\sum
_{j=0}^{c-5}a_j-b_1-a_{t+c-3}-a_{t+c-2}$ and the smallest possible
degree of a generator of $I_{c-1}$ is $\ell _c
-\sum_{j=t-1}^{t+c-3}a_j$. Hence  $\Ext
^1_R(I_{D_{c-2}},I_{c-1})_0=0$ if $\ell _c-\sum
_{j=0}^{c-5}a_j-b_1-a_{t+c-3}-a_{t+c-2}<\ell _c
-\sum_{j=t-1}^{t+c-3}a_j$ or, equivalently, $a_{t+c-2}>\sum
_{j=t-1}^{t+c-4}a_{j}-\sum _{j=0}^{c-5} a_j-b_1$ which is our
assumption  $(j_c).$

Similarly we get $\Ext ^1_R(I_{D_{i-1}},I_i)_0=0$ if $(j_{i+1})$
holds. Since we by Remark \ref{dep} and the hypothesis
$a_{i-min(3,t)}\ge b_{i}$ for $min(3,t)\le i \le t$  know that a
general enough $C\in W(\underline{b};\underline{a})$ satisfies
$\depth _{I(Z_{i})}D_i\ge 4$ for $2\le i \le c-2$, we conclude by
(\ref{510}).  For the dimension formula we use Remark \ref{debil}
(1).
\end{proof}

\vskip 2mm

Since Corollary \ref{component6} does not apply to $n=1$ and $3\le
c \le 4$, we include one more result to cover these cases. For
$c=3$, the result is known (\cite{KMMNP}, Corollary 10.11).

\begin{corollary}\label{c3c4} Let $W(\underline{b};\underline{a})$ be
the locus of good determinantal schemes in $\PP^{n+c}$ of
dimension $n\ge 1$. If either
\begin{itemize}
\item[(1)] $c=3$, $a_{i-min(2,t)}\ge b_{i}$ for $min(2,t)\le i \le t$
and $a_{t+1}>a_{t-1}+a_t-b_1$, or

\item[(2)] $c=4$, $a_{i-min(3,t)}\ge b_{i}$ for $min(3,t)\le i \le t$  and
$a_{t+2}>a_{t-1}+a_t-b_1$,
\end{itemize}
then $\overline{W(\underline{b};\underline{a})}$ is a generically
smooth, irreducible component of $\Hi ^p(\PP^{n+c})$ of dimension
$\lambda _c+K_3+...+K_c$.
\end{corollary}
\begin{proof} Let $c=3$. To see the vanishing of
$\Ext^1_{D_2}(I_{D_{2}}/I^2_{D_{2}},I_2)_0$ of Theorem
\ref{main5}, it suffices to prove $\Ext^1_R(I_{D_{2}},I_2)_0=0$.
As in the proof of Corollary \ref{component6},  we find the
minimal degree of relations of $I_{D_{2}}$ to be $\ell _2-b_1$,
and we get the vanishing of the $\Ext ^1_R$-group above by
assuming $\ell _2-b_1=\ell _3 -a_{t+1}-b_1<\ell _3-\sum_{j=t-1}^t
a_j$, i.e. $a_{t+1}>a_{t-1}+a_t-b_1$.

If $c=4$, it suffices to prove that $\Ext^1_R(I_{D_{2}},I_3)_0=0$
by the argument of (\ref{510}). Indeed, (\ref{53}) vanishes and we
know that $\depth _{I(Z_{2})}D_2\ge 4$ implies an injection
$$\Ext^1_{D_3}(I_{D_{3}}/I^2_{D_{3}},I_3)_0\hookrightarrow
\Ext^1_{D_2}(I_{D_{2}}/I^2_{D_{2}},I_3)_0$$ by (\ref{48new}) and
that the latter $\Ext ^1$-group vanishes if
$\Ext^1_R(I_{D_{2}},I_3)_0=0$. Now exactly as in the first part of
the proof of Corollary \ref{component6}, we have
$\Ext^1_R(I_{D_{2}},I_3)_0=0$ provided $$ (j_4): \quad
a_{t+2}>a_{t-1}+a_t-b_1. $$ and we conclude by Theorem
\ref{main5}. For  the dimension formulas, we use Remark
\ref{debil}(2) and Corollary \ref{cod3dim0}.
\end{proof}

\vskip 4mm Corollaries \ref{cod6} and \ref{component6} and Remark
\ref{debil}(1) can be improved a little bit if we increase $\depth
_{I(Z_i)}D_i$. In fact,  by Remark \ref{dep} we know that under
the assumption $a_{i-min(3,t)}\ge b_{i}$ for $min(3,t)\le i \le
t$, we can suppose $\depth _{I(Z_i)}D_i\ge 5$ for $i\ge 3$,
letting $Z_i=Sing(X_i)$. Since $Z_i\subset Z_{i-1}$, if we suppose

\begin{equation}\label{511}
\depth_{I(Z_{i-2})}D_{i-2}\ge 5
\end{equation}
we get $\depth_{I(Z_{i-2})}D_{i-1}\ge 4$, and hence $
\depth_{I(Z_{i-1})}D_{i-1}\ge 4$, as well as $
\depth_{I(Z_{i-2})}D_{i}\ge 3$. As in (\ref{new}), we see that

$$\Ext  ^1_{D_{i-2}}(I_{i-2},I_i)\cong H^1_*(U_{i-2},\cH
om_{\cO_{X_{i-2}}}(\cI_{i-2},\cI_i))\cong $$ $$H^1_*(U_{i-2},\cH
om_{\cO_{X_{i}}}(\cI_{i-2}\otimes \cI ^*_i,\cO_{X_i}))\cong
H^1_*(U_{i-2},\cO_{X_{i}}(a_{t+i-3}-a_{t+i-2}))=0. $$ Arguing as
in (\ref{47new}) and  combining with (\ref{48new}) we get
injections
\begin{equation}\label{512}
\Ext ^1_{D_i}(I_{D_i}/I^2_{D_i},I_i)\hookrightarrow \Ext
^1_{D_{i-1}}(I_{D_{i-1}}/I^2_{D_{i-1}},I_i)\hookrightarrow \Ext
^1_{D_{i-2}}(I_{D_{i-2}}/I^2_{D_{i-2}},I_i).
\end{equation}

In particular, if $\Ext ^1_R(I_{D_{i-2}},I_i)_0=0$ and if
(\ref{511}) hold, then $\Ext ^1_{D_i}(I_{D_i}/I^2_{D_i},I_i)_0=0$.
Now looking to the proof of Corollary \ref{component6}, we easily
see that $\Ext ^1_R(I_{D_{c-3}},I_{c-1})_0=0$ provided
\begin{equation}\label{513}
a_{t+c-2}>\sum_{j=t-1}^{t+c-5}a_j-\sum_{j=0}^{c-6}a_j-b_1.
\end{equation}

\noindent Hence if $(j_{c-1})$ holds, then (\ref{513})  holds
because $a_{t+c-2}\ge a_{t+c-3}$, and it is superfluous to assume
$(j_c)$ in Corollary  \ref{component6}. The argument requires
$\depth _{I(Z_{c-3})}D_{c-3}\ge 5$, i.e. $c\ge 6$, and arguing
slightly more general, we see that for $6\le i \le c$, then
$(j_i)$ is superfluous provided $(j_{i-1})$ holds.

In particular the conclusion of Corollary \ref{component6} holds
if $(j_i)$ holds for any {\em odd} number $i$ such that $5\le i
\le c$.

\begin{remark} (1) In Corollary \ref{component6}, the assumption
$(j_i)$ is superfluous if $(j_{i-1})$ holds, $ 6\le i \le c$.

 (2)
Increasing $\depth _{I(Z_i)}D_i$ even more (say by assuming
$a_{i-min(4,t)}\ge b_{i}$ for $min(4,t)\le i \le t$, cf. Remark
\ref{dep}), we can weaken $(j_k)$, resp. $(i_k)$, conditions of
Corollary \ref{component6}, resp. Corollary \ref{cod6}, further.
\end{remark}

\section{Conjecture}

\vskip 4mm We end the paper with a Conjecture raised by this
paper. In fact Theorem \ref{upperbound}, Proposition \ref{bound2},
and Corollaries \ref{cod3}, \ref{cod4}, \ref{cod5}, \ref{cod6} and
\ref{cod3dim0} suggest - and prove in many cases - the following
conjecture:

\begin{conjecture}
Given integers $a_0\le a_1\le ... \le a_{t+c-2}$ and $b_1\le
...\le b_t$, we set $\ell :=\sum_{j=0}^{t+c-2}a_j-\sum_{i=1}^tb_i$
and $h_i:= 2a_{t+1+i}+a_{t+2+i}+\cdots +a_{t+c-2}-\ell +n+c$, for
$i=0,1,...,c-3$. Assume $a_{i-min([c/2]+1,t)}\ge b_{i}$ for
$min([c/2]+1,t)\le i \le t$. Then we have

 \[
 \dim W(\underline{b};\underline{a}) = \sum_{i,j}
\binom{a_i-b_j+n+c}{n+c}- \sum _{i,j} \binom{a_i-a_j+n+c}{n+c}-
\]
\[
 \sum _{i,j} \binom{b_i-b_j+n+c}{n+c} + \sum _{j,i}
\binom{b_j-a_i+n+c}{n+c}+\binom{h_0}{n+c} + 1+ \]
\[ \sum_{i=1}^{c-3}\left (\sum
_{r+s=i} \sum _{0\le i_1< ...< i_{r}\le t+i \atop 1\le
j_1\le...\le j_s \le t } (-1)^{i-r} \binom{h_i+a_{i_1}+\cdots
+a_{i_r}+b_{j_1}\cdots +b_{j_s} }{n+c}\right  ).
\]

\end{conjecture}

In particular, we would like to know if the Conjecture 6.1 is at
least  true when  the entries of $\cA $ all have the same degree.
More precisely,

\begin{conjecture} Let $W(\underline{0};\underline{d})$ be the
locus of good determinantal schemes in $\PP^{n+c}$ of codimension
$c$ given by the maximal minors of a $t\times (t+c-1)$ matrix with
entries homogeneous forms of degree $d$. Then, $$\dim
W(\underline{0};\underline{d}) = t(t+c-1)\binom{d+n+c}{n+c}-
t^2-(t+c-1)^2+1.$$
\end{conjecture}

Finally, since \cite{KMMNP}, \S 10 has served as standard
reference, we take the opportunity to mention an inaccuracy in
\cite{KMMNP}, (10.12) and correct it.

\begin{remark} We propose to substitute the hypothesis  \cite{KMMNP},
(10.12):

\vskip 2mm {\em Given $C\subset \PP^{n+c}$  a good standard
determinantal scheme of dimension $n$, there always exists a flag
$C=X_c\subset X_{c-1}\subset ...\subset X_2\subset \PP^{n+c}$ such
that for each $i<c$, the closed embedding $X_i \hookrightarrow
\PP^{n+c}$ and $X_{i+1}\hookrightarrow X_i$ are local complete
intersection (l.c.i.) outside some set $Z_i$ of codimension 1 in
$X_{i+1}$ ($\depth _{Z_i}\cO _{X_{i+1}}\ge 1$).}

\vskip 2mm
 \noindent by (*):

\vskip 2mm
 {\em Given $C\subset \PP^{n+c}$  a good standard determinantal
scheme of dimension $n$, we will assume that there exists a flag
$C=X_c\subset X_{c-1}\subset ...\subset X_2\subset \PP^{n+c}$ such
that for each $i<c$, the closed embedding $X_{i+1}\hookrightarrow
X_i$ is l.c.i. outside some set $Z_i$ of codimension 2 in
$X_{i+1}$ ($\depth _{Z_i}\cO _{X_{i+1}}\ge 2$). Moreover, we
suppose $X_2\hookrightarrow \PP^{n+c}$ is a l.c.i. in codimension
$\le 1$. }

\vskip 2mm The reason of increasing the depth related to
$X_{i+1}\hookrightarrow X_i$ by 1 is that the exactness of
\cite{KMMNP},  (10.15) in the proof of \cite{KMMNP}, Proposition
10.12 is straightforward to see if (*) holds (by e.g. (\ref{NM}))
while it is doubtful with the original hypothesis. Hence in
\cite{KMMNP}, Proposition 10.12, Theorem 10.13, Remarks 10.6 and
10.14, Example 10.16, Corollary 10.17 ($n=1$) and Example 10.18
($n=1$) we should suppose (*) instead of \cite{KMMNP}, (10.12). So
in \cite{KMMNP} Corollary 10.17 ($n=1$) and Example 10.18 ($n=1$)
we need to suppose $C\subset S$ to be Cartier instead of
generically Cartier and $S$ to be $G_1$, while \cite{KMMNP},
Corollary 10.15 need no change because $C\subset S$ is supposed to
be Cartier outside a sufficiently small $Z$. In \cite{KMMNP},
Corollary 10.15 we {\em may} replace \cite{KMMNP}, (10.12) by (*)
and hence by "$S$ is $G_1$", and hence we point out that
\cite{KMMNP}, Propositions 10.7 and 1.12 are valid as stated.

\vskip 4mm Now we consider the 0-dimensional case. Looking closer
to the proof of \cite{KMMNP}, Proposition 10.12, we need the
graded version of \cite{KMMNP} (10.15) to be exact in degree zero.
Therefore still assuming (*) the proof is only complete for a flag
$R\twoheadrightarrow D_2\twoheadrightarrow ...\twoheadrightarrow
D_c=A$ satisfying $\dim A\ge 2$ or $\dim A=1$ and
$M_c(a_{t+c-2})_0\cong (N_{D_c/D_{c-1}})_0$, cf. (\ref{NM}). Hence
in \cite{KMMNP}, Theorem 10.13 the $H$-unobstructedness (when
$\dim C=0$) does not necessarily follow from \cite{KMMNP},
Proposition 10.12. Fortunately, \cite{KMMNP}, Theorem 9.6, makes
explicit an assumption which implies the $H$-unobstructedness of
$C$ ( (iii) below), and we can weaken (iii) further to $\Ext
^1_A(I_c/I_c^2,A)_0=0$ by the proof of \cite{KMMNP}, Theorem 9.6
because  $H^2(D_{c-1},A,A)\cong \Ext ^1_A(I_c/I_c^2,A)$. Summing
up, in the $\dim C=0$ case of \cite{KMMNP}, Theorem 10.13 and
Corollary 10.17, if we assume (*) instead of \cite{KMMNP},
(10.12), and in addition at least one of the following conditions

\begin{itemize}
\item[(i)] $C$ is $H$-unobstructed (e.g.$\Ext ^1_A(I_c/I_c^2,A)_0=0$),
\item[(ii)]$M_c(a_{t+c-2})_0\cong  (N_{A/D_{c-1}})_0$ (i.e. (\ref{DiMi})
extends to a short exact sequence in degree zero for $i=c-1$),
\item[(iii)] $\Hom _R(I_c,H^2_{\goth{m}}(I_A))_0=0$,
\end{itemize}

\noindent then the conclusions  hold. The only example in which
$\dim C=0$ is \cite{KMMNP}, Example 10.18. In this example (iii)
holds because the degree of all minimal generators of $I_C$ are
$m-t-1$ and $H^2_{\goth{m}}(I_A)_{\nu} \cong H^1(\cI_C(\nu))=0$
for $\nu
>m+t-3$, cf. \cite{KMMNP}, Example 7.5 for details.
Similarly, in the situation of \cite{KMMNP}, Corollary 10.17,
using the order $b_1\ge ... \ge b_t$ and $a_0\ge a_1\ge ... \ge
a_{t+1}$ appearing in \cite{KMMNP}, Corollary 10.17 we see that
(iii) holds provided $a_1+a_2\le 3+  2b_t$, a rather strong
condition. The condition (i)
might however be week.
\end{remark}

\vskip 2mm So  while the substitution of \cite{KMMNP} (10.12) by
(*) in the case $\dim C\ge 1$ is relatively innocent (i.e.
$X_{i+1}\hookrightarrow X_i$ has to be Cartier in codimension $\le
1$ instead of generically Cartier), the 0-dimensional case leads
to an extra assumption. Therefore, it is probably a more natural
approach in the 0-dimensional case to just find $\dim
W(\underline{b}; \underline{a})$ without  trying to prove that $
W(\underline{b}; \underline{a})$ is an irreducible component of
GradAlg and if $\dim C\ge 1$, it is natural both to find $\dim
W(\underline{b}; \underline{a})$ and to see that $
W(\underline{b}; \underline{a})$ is a generically smooth,
irreducible component of $\Hi ^p(\PP^{n+c})$ as it has been the
strategy of this paper.

\end{document}